\documentstyle[12pt]{article}
\title{Optimal control on distributions}
\author{Constantin Udri\c ste}
\begin{document}
\maketitle

\newcommand{\area}{\hbox{area}}
\newcommand{\ome}{\Omega}

\newcommand{\omm}{\Omega}
\newcommand{\grad}{\hbox{$\:$grad \,}}
\newcommand{\divv}{\hbox{$\:$div \,}}
\newcommand{\rot}{\hbox{$\:$rot \,}}
\newcommand{\im}{\hbox{$\:$Im }}
\newcommand{\re}{\hbox{$\:$Re }}
\newcommand{\sh}{\hbox{$\:$sh\, }}
\newcommand{\ch}{\hbox{$\:$ch\,}}
\newcommand{\om}{\omega}
\newcommand{\br}{\hbox{\bf R}}
\newcommand{\bc}{\hbox{\bf C}}
\newcommand{\bz}{\hbox{\bf Z}}
\newcommand{\pp}{\prime}
\newcommand{\ty}{\infty}
\newcommand{\di}{\displaystyle}
\newcommand{\va}{\varphi}
\newcommand{\si}{\sigma}
\newcommand{\ga}{\gamma}
\newcommand{\gaa}{\Gamma}
\newcommand{\na}{\nabla}
\newcommand{\te}{\theta}
\newcommand{\ld}{\ldots}
\newcommand{\ov}{\over}
\newcommand{\ri}{\Rightarrow}
\newcommand{\noa}{\noalign{\medskip}}
\newcommand{\la}{\lambda}
\newcommand{\su}{\subset}
\newcommand{\qu}{\quad}
\newcommand{\fo}{\forall}
\newcommand{\al}{\alpha}
\newcommand{\lll}{\Leftrightarrow}
\newcommand{\be}{\beta}
\newcommand{\ep}{\varepsilon}
\newcommand{\pa}{\partial}
\newcommand{\ti}{\times}
\newcommand{\dd}{\Delta}
\newcommand{\sgg}{\Sigma}
\newcommand{\de}{\delta}


\newcommand\C{{\,I\!\!\!\!\!\: C}}
\newcommand{\ds}{\displaystyle}
\newcommand{\ol}{\overline}
\newcommand{\w}{\widetilde}
\newcommand{\ul}{\underline}
\newcommand{\s}{\stackrel}
\newcommand{\p}{\prime}
\def\d{\displaystyle\frac}
\def\pa{\partial}
\def\wi{\widetilde}
\def\ov{\over}
\def\ri{\rightarrow}
\def\ba{\begin{array}}
\def\ea{\end{array}}
\def\Ker{{\rm{Ker}}\,}
\def\var{\varepsilon}
\def\fii{\varphi}
\def\Lam{\Lambda}
\def\lam{\lambda}
\def\ii{\^\i}
\def\dx{\dot x}
\def\ddx{\ddot x}
\def\qu{\quad}
\def\al{\alpha}
\def\sul{\displaystyle\sum\limits}
\def\Lr{\Leftrightarrow}
\def\we{\wedge}
\def\I{\^I}
\def\a{\^a}
\def\lri{\longrightarrow}
\def\Lri{\Longrightarrow}
\def\Llr{\Longleftrightarrow}
\def\llr{\longleftrightarrow}
\def\br{R}
\def\om{\omega}
\def\qu{\quad}
\def\ld{\ldots}
\def\fo{\forall}
\def\di{\displaystyle}
\def\al{\alpha}
\def\ep{\epsilon}
\def\su{\subset}
\def\gaa{\Gamma}
\def\ty{\infty}
\def\va{\varphi}
\def\noa{\noalign}
\def\ti{\times}
\def\pp{\prime}

\setcounter{page}{1}
\textheight=19cm \textwidth=13cm
\oddsidemargin=16mm \evensidemargin=16mm
\date{}


{ \footnotesize \it This paper studies (single-time and multitime) optimal control problems
on a nonholonomic manifold (described either by the kernel of a Gibbs-Pfaff form or by
the span of appropriate vector fields). For both descriptions we analyse: infinitesimal deformations
and adjointness, single-time optimal control problems, multitime optimal control
problem of maximizing a multiple integral functional, multitime optimal control problem of
maximizing a curvilinear integral functional, Curvilinear functionals depending on curves,
optimization of mechanical work on Riemannian manifolds.
Also we prove that a nonholonomic system can be always controlled by
uni-temporal or bi-temporal bang-bang controls.}

{\bf Mathematics Subject Classification 2010}: 49J15, 49J20, 93C15, 93C20.

{\bf Keywords:} nonholonomic manifold,  single-time optimal control, multitime optimal control, bang-bang controls.

\section{Optimal control on a distribution \\described by a Pfaff equation}

A {\it generalized distribution} $\Delta: x \to \Delta_x$, or Stefan-Sussmann distribution, is similar to a distribution,
but the subspaces are not required to be all of the same dimension. The definition
requires that the subspaces $\Delta_x$ are locally spanned by a set of vector fields,
but these will no longer be everywhere linearly independent. It is not hard to see that the dimension of the distribution $\Delta$ is
lower semicontinous, so that at special points the dimension is lower than at nearby points.
One class of examples is provided by a non-free action of a Lie group on a manifold,
the vector fields in question being the infinitesimal generators of the group
action (a free action gives rise to a genuine distribution).
Another examples arise in dynamical systems, where the set of vector fields in
the definition is the set of vector fields that commute with a given one.
There are also examples and applications in Control theory, where the
generalized distribution represents infinitesimal constraints of the system.

see ControlJakubczyk, pag 146

{\bf Lemma} {\it If the variational system (treated as linear system without constraints on the control) is controllable, then
the original system is strongly accessible}.

Nonholonomic path planning represents a fusion of some of the
newest ideas in control theory, classical mechanics, and differential geometry
with some of the most challenging practical problems in robot
motion planning. Furthermore, the class of systems to which the theory
is relevant is broad: mobile robots, space-based robots, multifingered
hands, and even such systems as a one-legged hopping robot. The
techniques presented here indicate one possible method for generating
efficient and computable trajectories for some of these nonholonomic
systcms in the absence of obstacles.

The delay differential equations (DDEs) are a type of differential equation
in which the derivative of the unknown function at a certain time is given in terms of the values of the function at previous times.
For example
$$\frac{d}{dt}x(t) = f(x(t), x(t-\tau)).$$
A delay Pfaff equation means
$$a_i(x(t),x(t - \tau))dx^i(t) = 0.$$

\subsection{Infinitesimal deformations and \\adjointness on distributions}

Let $D$ be a nonholonomic distribution on $R^n$ described by a Pfaff equation
$$a_i(x)dx^i = 0.\leqno(1)$$

Let $x(t),\,\,t\in I = [t_0, t_1] \subset R$, be an integral curve of the distribution $D$.
Let $x(t;\epsilon),\,\epsilon \in [0,\de)$ be a differentiable variation of $x(t)$, i.e.,
$$A_i(x(t;\epsilon); \epsilon)dx^i(t;\epsilon) = 0,\, A_i(x;0) = a_i(x), \,x(t;0) = x(t).$$
The variation $(t,\ep) \to x(t;\epsilon),\, t \in I,\,\epsilon \in [0,\de)$ is a surface. It is an integral surface
only if the distribution $D$ admits integral surfaces.
Taking the partial derivative with respect to $\epsilon$
and denoting $y^i(t) =\di \frac{\pa x^i}{\pa \epsilon}(t;0)$, we find the {\it single-time (Pfaff) infinitesimal deformation equation}
$$\left(\frac{\pa a_i}{\pa x^j}(x)y^j(t) + b_i(x)\right)dx^i + a_i(x) dy^i = 0,\, b_i (x)= \frac{\pa A_i}{\pa \epsilon}(x;0)\leqno(2)$$
around a solution $x(t)$ of the Pfaff equation (1).
The {\it single-time adjoint Pfaff system} is
$$d(p(t)a_j(x(t))) =\frac{\pa (p(t) a_i(x)dx^i)}{\pa x^j},\leqno(3)$$
whose solution $p(t)$ is called the {\it costate function}. The foregoing Pfaff equations (2) and (3) are {\it adjoint} ({\it dual})
in the following sense: if $y$ is a solution of the infinitesimal deformation (Pfaff) equation (2), then the function
$p(t)a_i(x)y^i(t)$ verifies the Pfaff equation $d(p(t)\,a_i(x) y^i(t)) + p(t)\,b_i(x)dx^i = 0$.

Let $x(t),\,\, t \in \Omega_{t_0t_1}\subset R^m_+$ be a maximal, $m$-dimensional, $m \leq n-2$,
integral submanifold of the distribution $D$.
We fix $k = 1,...,n- m -1$. Let $\epsilon = (\epsilon^A)\in [0,\de)^k, \, A = 1,...,k,$
and let $x(t;\epsilon)$ be a differentiable variation of $x(t)$, i.e.,
$$A_i(x(t;\epsilon);\epsilon)dx^i(t;\epsilon) = 0,\, A_i(x,0) = a_i(x),\, x(t;0) = x(t).$$
The variation $(t,\ep) \to x(t;\epsilon),\, t \in \Omega_{t_0t_1},\,\epsilon \in [0,\de)^k$ is an $(m+k)$-dimensional manifold,
but not an integral submanifold.

Taking the partial derivative with respect to $\epsilon^A$
and denoting $y^i_A(t) = \frac{\pa x^i}{\pa \epsilon^A}(t;0)$, we find the {\it multitime infinitesimal deformation (Pfaff) system}
$$\left(\frac{\pa a_i}{\pa x^j}(x)y^j_A(t) + b_{iA}(x)\right)dx^i + a_i(x) dy^i_A = 0,\,b_{iA} (x)= \frac{\pa A_i}{\pa \epsilon^A}(x;0)\leqno(4)$$
around a solution $x(t)$ of the Pfaff equation (1). The {\it multitime adjoint Pfaff system} is
$$d(p^A(t)a_j(x)) = \frac{\pa (p^A(t)a_i(x)dx^i)}{\pa x^j},\leqno(5)$$
whose solution $p(t) = (p^A(t))$ is called the {\it costate vector}. The foregoing Pfaff equations (4) and (5) are {\it adjoint} ({\it dual})
in the following sense: if $y^i_A$ is a solution of the infinitesimal deformation (Pfaff) system (4), then the function
$a_i(x)p^A(t)\,y^i_A(t)$ verify the Pfaff equation $d(a_i(x) p^A(t)\,y^i_A(t)) + p^A(t)\,b_{iA}(x)dx^i= 0$.

\subsection{Evolution of a distribution}

Let $D$ be a nonholonomic distribution on $R^n$ described by a Pfaff equation (1).
Let $x(t),\,\,t\in I = [t_0, t_1] \subset R$, be an integral curve of the distribution $D$.
Let $x(t;\epsilon),\, \ep \in [0,\de)$ be a differentiable variation of $x(t)$. Suppose that $(x(t;\epsilon);\epsilon)$
is an integral surface of a Pfaff equation in $R^{n+1}$, i.e.,
$$A_i(x(t;\epsilon), \epsilon)dx^i(t;\epsilon) = B(x(t;\epsilon), \epsilon)d\epsilon,$$
$$A_i(x,0) = a_i(x), \,B(x,\epsilon) = A_i(x,\epsilon)\frac{\pa x^i}{\pa \epsilon},\, B(x,0) = 0,\,x(t;0) = x(t).$$
Taking the partial derivative with respect to $\epsilon$, we find
$$\left(\frac{\pa A_i}{\pa x^j} \frac{\pa x^j}{\pa \epsilon} + \frac{\pa A_i}{\pa \epsilon}\right) dx^i + A_i\, d\, \frac{\pa x^i}{\pa \epsilon}
= \left(\frac{\pa B}{\pa x^j} \frac{\pa x^j}{\pa \epsilon} + \frac{\pa B}{\pa \epsilon}\right) d\epsilon.$$
If we accept an evolution after the direction of the vector field $X^j$, i.e.,
$\di\frac{\pa x^j}{\pa \epsilon} = \al\, X^j$, then we find the PDE system
$$\left(\frac{\pa A_i}{\pa x^j} X^j + \frac{\pa A_i}{\pa \epsilon}\right) dx^i + A_i\, dX^i
= \left(\frac{\pa B}{\pa x^j} X^j + \frac{\pa B}{\pa \epsilon}\right) d\epsilon$$
with unknowns $A_i$, fixed by initial conditions. For $\epsilon = 0$, we rediscover the system in variations,
with the condition $a_i(x)y^i(t)=0, \, y^i(t) = \frac{\pa x^i}{\pa \ep}(t,0)$.

\subsection{Single-time optimal control problems on \\a distribution}

Let $D$ be a distribution on $R^n$ described by a controlled Pfaff equation
$$a_i(x, u)dx^i = 0,\, x = (x^i)\in R^n,\, u = (u^a)\in R^k$$
and let $x(t),\,\,t\in I = [t_0, t_1],$ be an integral curve of the distribution $D$.

{\it A  single-time optimal control problem consists of  maximizing the functional
$$I(u(\cdot)) = \int_{t_0}^{t_1}L(t,x(t),u(t)) dt +g(x(t_1))\leqno(6)$$
subject to
$$a_i(x(t), u(t))dx^i(t) = 0,\,\hbox{a.e.}\,\,\, t\in I = [t_0, t_1],\, x(t_0) = x_0.\leqno(7)$$}
It is supposed that $L: I \times A\times U \to R$ is a $C^2$ function, $a_i: A\times U \to R, \, i = 1,..., n$
are $C^2$ functions and $g(x(t_1))$ is a $C^1$ function.
Ingredients: $A$ is a bounded and closed subset of $R^n$ which contains each trajectory $x(t),\, t\in I$
of controlled system, and $x_0$ and $x_1$ are the initial and final states of the trajectory $x(t)$.
The values of the control functions belong to a set $U\subset R^k$, bounded and closed.

Let us find the first order necessary conditions for an optimal pair $(x,u)$. We fix the control $u(t)$ and we variate
the state $x(t)$ into $x(t,\ep)$. We obtain the {\it single-time infinitesimal deformation (Pfaff) equation}
$$\frac{\pa a_i}{\pa x^j}(x,u)y^j dx^i + a_i(x,u) dy^i = 0.$$
of the nonholonomic constraint $a_i(x(t), u(t))dx^i(t) = 0$.
It follows the {\it single-time adjoint Pfaff equation}
$$d(p(t) a_j(x(t),u(t))) = \frac{\pa (p(t) a_i(x,u)dx^i)}{\pa x^j},$$
whose solution $p(t)$ is called the {\it costate function}.
Here, the symbol $d$ in the left hand member of the adjoint equation
means the differentiation with respect to $p$ and $x$.

Using the Lagrangian $1$-form
$${\cal L} = L(t,x(t),u(t)) dt,$$
we build the {\it Hamiltonian $1$-form}
$${\cal H} = {\cal L}(t,x(t),u(t)) + p(t) a_i(x(t), u(t))dx^i(t).$$

{\bf Theorem} ({\bf Single-time maximum principle})
{\it Suppose that the problem of maximizing the functional (6) constrained by (7) has
an interior optimal solution $\hat u(t)$, which determines
the optimal evolution $x(t)$. Then there exists a
costate function $(p(t))$ such that
$$
\frac{\partial {\mathcal H}}{\partial p}= a_i(x(t), u(t))dx^i(t) = 0, \leqno(8)
$$
the function $(p(t))$ is the unique solution of the following Pfaff system (adjoint system)
$$
d(p(t)a_j(x(t),u(t))) = \frac{\pa {\cal H}}{\pa x^j}\leqno(9)
$$
and satisfies the critical point conditions
$$
{\mathcal H}_{u^a}\left( t,x(t),u(t),p(t)\right) =0,\,a = \overline{1,k}.\leqno(10)
$$}

{\bf Proof} We use the Hamiltonian $1$-form ${\cal H}$.
The solutions of the foregoing problem are among the solutions of the free maximization problem of the curvilinear integral functional
$$J(u(\cdot)) = \int_{\tilde\Gamma}{\cal H}(t,x(t),u(t),p(T)) + g(x(t_1)),$$
where $\tilde \Gamma = ([t_0,t_1], x([t_0,t_1]))=\{(t,x(t))|t \in [t_0,t_1]\}\subset R_+\times R^n$.

Suppose that there exists a continuous control $\hat u(t),\, t \in I = [t_0,t_1]$, with $\hat u(t)\in \hbox{Int}\,{U}$,
and an integral curve $x(t)$ which are optimal in the previous problem. Now consider a control variation
$u(t,\epsilon)=\hat u(t) + \epsilon h(t),\, t \in I = [t_0,t_1]$, where $h$ is an arbitrary continuous vector function,
and a state variation $x(t,\ep),\, t \in I = [t_0,t_1]$, related by
$$a_i(x(t;\epsilon), u(t; \epsilon))dx^i(t;\epsilon) = 0,\, x(t;0) = x(t),\, x(t_0,0)=x_0.$$
Since $\hat u(t)\in \hbox{Int}\,{\cal U}$ and any continuous function over a compact set
$I$ is bounded, there exists $\epsilon_h>0$ such
that $u(t,\epsilon)=\hat u(t) + \epsilon h(t)\in \hbox{Int}\,{U},\,\,\forall |\epsilon|<\epsilon_h$.
This $\epsilon$ is used in our variational arguments.

For $|\epsilon|<\epsilon_h$, we define the function
$$J(\epsilon)= \int_{\tilde\Gamma(\ep)}{\cal H}(t,x(t,\epsilon),u(t,\epsilon),p(t)) + g(x(t_1,\ep))$$
$$= \int_{\tilde\Gamma(\ep)}{\cal L}(t,x(t,\epsilon),u(t,\epsilon)) + p(t) a_i(x(t,\epsilon), u(t,\epsilon)) dx^i(t,\epsilon)+ g(x(t_1,\ep)).$$
Differentiating with respect to $\epsilon$, it follows
$$
J^{\prime}(\epsilon)= \int_{\tilde\Gamma(\ep)}{\cal L}_{x^j}(t,x(t,\epsilon),u(t,\epsilon))x^j_{\epsilon}(t,\epsilon) + \frac{\pa g}{\pa x^j}(x(t_1,\ep))x^j_\ep(t_1,\ep)$$
$$ + \int_{\tilde\Gamma(\ep)} p(t) \frac{\pa a_i}{\pa x^j}(x(t,\epsilon), u(t,\epsilon))x^j_{\epsilon}(t,\epsilon) dx^i(t,\epsilon)$$
$$ +  \int_{\tilde\Gamma(\ep)} p(t) a_i(x(t,\epsilon), u(t,\epsilon)) dx^i_{\epsilon}(t,\epsilon)$$
$$ + \int_{\tilde\Gamma(\ep)}{\cal L}_{u^a}(t,x(t,\epsilon),u(t,\epsilon)h^a(t) $$
$$+ \int_{\tilde\Gamma(\ep)} p(t) \frac{\pa a_i}{\pa u^a}(x(t,\epsilon), u(t,\epsilon))h^a(t) dx^i(t,\epsilon)).$$
Evaluating at $\epsilon = 0$, we find $\tilde\Gamma(0) = \tilde\Gamma$ and
$$J^{\prime}(0)= \int_{\tilde\Gamma}\left({\cal L}_{x^j}(t,x(t),\hat u(t))+ p(t) \frac{\pa a_i}{\pa x^j}(x(t), \hat u(t))dx^i(t)\right)x^j_{\epsilon}(t,0)$$
$$ +  \int_{\tilde\Gamma} p(t) a_i(x(t), \hat u(t)) dx^i_{\epsilon}(t,0) + \frac{\pa g}{\pa x^j}(x(t_1,0))x^j_\ep(t_1,0)$$
$$ + \int_{\tilde\Gamma}\left({\cal L}_{u^a}(t,x(t),\hat u(t)) + p(t) \frac{\pa a_i}{\pa u^a}(x(t), \hat u(t))dx^i(t)\right) h^a(t),$$
where $x(t)$ is the curve of the state variable corresponding to the optimal control $\hat u(t)$. The integral from the middle can be written
$$\int_{\tilde\Gamma} p(t) a_j(x(t), \hat u(t)) dx^j_{\epsilon}(t,0) = p(t) a_jx^j_{\epsilon}|_{\pa{\tilde\Gamma}} - \int_{\tilde\Gamma}d(p(t)a_j(x(t), \hat u(t)))x^j_{\epsilon}(t,0),$$
where the symbol $d$ in the last integral means the differentiation with respect to $p$ and $x$.
We find $J^{\prime}(0)$ as
$$\int_{\tilde\Gamma}\left({\cal L}_{x^j}(t,x(t),\hat u(t)) + p(t)\frac{\pa a_i}{\pa x^j}(x(t), \hat u(t))dx^i(t) - d(p(t)a_j(x(t), \hat u(t))\right)x^j_{\epsilon}(t,0)$$
$$+ p(t) a_j(x(t), \hat u(t))x^j_{\epsilon}(t,0)|_{t_0}^{t_1} + \frac{\pa g}{\pa x^j}(x(t_1,0))x^j_\ep(t_1,0)$$
$$ + \int_{\tilde\Gamma}\left({\cal L}_{u^a}(t,x(t),\hat u(t)) + p(t) \frac{\pa a_i}{\pa u^a}(x(t), \hat u(t))dx^i(t)\right) h^a(t).$$
Since $x(t_0)=x_0$, we have $x^j_{\epsilon}(t_0,0) = 0$.
We select the costate function $p(t)$ as solution of the adjoint Pfaff equation
$${\cal L}_{x^j}(t,x(t),\hat u(t)) + p(t)\frac{\pa a_i}{\pa x^j}(x(t), \hat u(t))dx^i(t) - d(p(t)a_j(x(t), \hat u(t)) = 0,$$
with the terminal condition $p(t_1) a_j(x(t_1), \hat u(t_1)) + \frac{\pa g}{\pa x^j}(x(t_1)) =0$.
On the other hand, we need $J^{\prime}(0)=0$ for all $h(t)=(h^a(t))$. Since
the variation $h$ is arbitrary, we get the following ({\it critical point condition})
$$\frac{\pa {\cal L}}{\pa {u^a}}(t,x(t),\hat u(t)) + p(t) \frac{\pa a_i}{\pa u^a}(x(t), \hat u(t))dx^i(t) = 0.$$

The foregoing equations (9) and (10) can be written
$$d(p(t)a_j(x(t),u(t))) = \frac{\pa {\cal H}}{\pa x^j}\left( t,x(t),u(t),p(t)\right),$$
$$\frac{\pa {\cal H}}{\pa {u^a}}\left( t,x(t),u(t),p(t)\right) = 0.$$

{\bf Example} Let us solve the problem
$$\max J(u(\cdot)) = - \frac{1}{2}\int_{t_0}^{t_1}(u^2(t) + z^2(t))dt$$
subject to (controlled Martinet distribution)
$$dz = \frac{1}{2}(y^2 + u)dx,\,x(t_0) = x_0, y(t_0) = y_0, z(t_0) = z_0.$$

Denote $\om = \frac{1}{2}(y^2 + u)dx -dz$. Since $$\om \wedge d\om = \frac{1}{2}dx \wedge d(y^2+u) \wedge dz$$
the distribution admits only integral curves with the parameter $t$. The Pfaff equation
$$dz = \frac{1}{2}(y^2 + u)dx,$$
is equivalent to a differential equation
$$\dot z(t) =  \frac{1}{2}(y(t)^2 + u(t))\dot x(t)$$
or to the ODE system
$$\dot x(t) = \varphi(t),\, \dot y(t) = \psi(t), \, \dot z(t)= \frac{1}{2}(y(t)^2 + u(t)) \varphi(t).$$

{\bf First variant} Using the Hamiltonian
$$H = - \frac{1}{2}(u^2(t) + z^2(t)) + p_1(t) \varphi(t) + p_2(t) \psi(t)) + p_3(t) \left(\frac{1}{2}(y(t)^2 + u(t))\varphi(t)\right)\!\!,$$
we find the adjoint ODEs
$$\dot p_1(t)  = - \frac{\pa H}{\pa x} = 0,\, \dot p_2(t)  = - \frac{\pa H}{\pa y}= - p_3 y \varphi,\, \dot p_3(t) = - \frac{\pa H}{\pa z} = z(t)$$
and the critical point condition
$$\frac{\pa H}{\pa u} = - u + \frac{1}{2}\,p_3 \varphi = 0.$$
We find the control $u = \frac{1}{2}\,p_3\varphi$ and $p_1(t) = c_1$. We need to find solutions for the system
$$\dot p_2(t) = - p_3(t) y(t) \dot x(t),\, \dot p_3(t) = z(t),\, \dot z(t) =  \frac{1}{2}\left(y(t)^2 + \frac{1}{2}p_3(t)\dot x(t)\right)\dot x(t).$$

{\bf Second variant} The Hamiltonian $1$-form is
$${\mathcal H} = - \frac{1}{2} (u^2(t) + z^2(t))dt + p(t) \left(\frac{1}{2}(y^2 + u)dx - dz\right).$$
The critical point condition
$$\frac{\pa {\mathcal H}}{\pa u} = - u dt+ \frac{1}{2}\,p\,dx = 0$$
gives the control $u = \frac{1}{2}\,p\dot x$.

Since $a_1 = \frac{1}{2}(y^2 + u)$, $a_2 = 0$, $a_3 = -1$, we find the adjoint Pfaff equations
$$d(p\frac{1}{2}(y^2 + u)) = \frac{\pa {\cal H}}{\pa x} = 0,\, 0 = \frac{\pa {\cal H}}{\pa y} = pydx, \,dp = - \frac{\pa {\cal H}}{\pa z} = zdt.$$
We need to solve the system
$$p(y^2 + u)= c, pydx = 0, dp = zdt,$$
$$udt = \frac{1}{2}\,p dx, dz = \frac{1}{2}(y^2 + u)dx.$$

On the other hand, the second variant ofers two explicit extremals:
(1) $y=0, pu=c_1, dp = zdt, udt = \frac{1}{2}pdx, dz = udx$, which does not satisfy the general
initial conditions; (2) $x=c_1, u = 0, z = c_3, p =c_3 t + c_4, y = \pm\sqrt{\frac{c}{c_3t+c_4}}$, depending upon four
arbitrary constants, which determine from initial conditions and terminal condition.

The first variant can be identified with the second variant via $p_1 = pa_1; p_2 = pa_2, p_3 = - pa_3$.

{\bf Third variant (Ionel Tevy)} We introduce two auxiliary controls $u_1, u_2$, changing the Pfaff equation into the controlled ODE system
$$\dot x(t) = u_1(t), \dot y(t) = u_2(t), \dot z (t)= \frac{1}{2}(y(t)^2 + u(t))u_1(t).$$
Then
$$H = -\frac{1}{2}(z^2 + u^2) + p_1u_1 + p_2u_2 + \frac{1}{2}p(y^2 + u)u_1.$$
We find the adjoint ODEs and the critical point conditions
$$\dot p_1 = - \frac{\pa H}{\pa x}=0, \dot p_2 = - \frac{\pa H}{\pa y} = -pyu_1, \dot p = z, $$
$$\frac{\pa H}{\pa u_1} = p_1 + \frac{1}{2}p(y^2 + u)=0, \frac{\pa H}{\pa u_2} = p_2 =0, \frac{\pa H}{\pa u} = - u +\frac{1}{2}pu_1=0.$$
It follows two extremals:
$$p_1(t) = c,\, p_2(t)=0,\, p(t) = -c_3 t - a,$$
$$u_1(t) = 0, \,u_2(t) = \dot y(t),\, u(t)=0,$$
$$x(t) = c_1, \,y(t) = \pm \sqrt{\frac{2c}{c_3t + a}}, \,z = - c_3;$$
$$p_1(t) = c, \,p_2(t) = 0, \,p(t) = \sqrt{\frac{(\al t +\be)^2 + 4c^2}{\al}},$$
$$u_1 = \frac{4\al c}{(\al t +\be)^2 + 4c^2}, \,u_2(t) =0, \,u(t) = - 2c \sqrt{\frac{\al}{(\al t +\be)^2 + 4c^2}},$$
$$x(t) = \int u_1(t)dt = 2 \,\hbox{atan}\frac{\al t + \be}{2c} + c_1, \,y(t) = 0,\, z(t) = \frac{\al t + \be}{\sqrt{(\al t +\be)^2 + 4c^2}}.$$

O varianta a rezolvarii Tevy se gaseste in Michele Pavon, Optimal control of nonholonomic systems, Proceedings of the 17th
International Symposium on Mathematical Theory of Networks and Systems, Kyoto, July 24-28, 2006; vezi NonhOptCon.pdf; citeaza lucrarile mele;

\subsection{Multitime optimal control problems on \\a distribution}

Let $D$ be a distribution on $R^n$ described by a controlled Pfaff equation
$$a_i(x, u)dx^i = 0,\, x = (x^i)\in R^n,\, u = (u^a)\in R^k$$
and let $x(t),\,\, t \in \Omega_{t_0t_1}\subset R^m$ be an $m$-dimensional
($m < n$, m maximal) integral submanifold of the distribution $D$.

Let us start with a

\subsubsection{ Multitime optimal control problem of \\maximizing a multiple integral functional}

{\it Find}
$$\max_{u(\cdot)}\,I(u(\cdot)) = \int_{\Omega_{t_0t_1}}L(t,x(t),u(t))\, \om + \int_{\pa{\Omega_{t_0t_1}}}g(x(t))d\si \leqno(11)$$
{\it subject to}
$$a_i(x(t), u(t))dx^i(t) = 0,\,\hbox{a.e.}\,\,\, t\in \Omega_{t_0t_1},\, x(t_0) = x_0.\leqno(12)$$
It is supposed that $L: \Omega_{t_0t_1} \times A\times U \to R$ is a $C^2$ function and $a_i: A\times U \to R, \, i = 1,..., n$
are $C^2$ functions. Ingredients: $\om = dt^1\cdots dt^m$ is the volume element, $A$ is a bounded and closed subset of $R^n$,
containing the images of the $m$-sheets $x(t),\,\ t\in \Omega_{t_0t_1}$ of controlled system,
and $x_0$ and $x_1$ are the initial and final states of an $m$-sheet $x(t)$.
The set in which the control functions takes their values is called as $U$, which is a bounded and
closed subset of $R^k$.

Let us find the first order necessary conditions for an optimal pair $(x,u)$.
We fix the control $u(t)$ and variate the state $x(t)$ into $x(t,\ep),\, \ep = (\ep^1,...,\ep^m)$.
We find the {\it multitime infinitesimal deformation (Pfaff) system}
$$\frac{\pa a_i}{\pa x^j}(x,u)y^j_\al(t) dx^i + a_i(x,u) dy^i_\al = 0$$
of the nonholonomic constraint $a_i(x(t), u(t))dx^i(t) = 0$. It follows the {\it multitime adjoint Pfaff system}
$$d(p^\al(t) a_j(x,u)) = \frac{\pa (p^\al(t)a_i(x,u)dx^i)}{\pa x^j},$$
whose solution $p(t) = (p^\al(t))$ is called the {\it costate vector}. Here,
the symbol $d$ in the left hand member of the adjoint equation
means the differentiation with respect to $p$ and $x$.

We use the Lagrangian $m$-form
$${\cal L} = L(t,x(t),u(t)) \om.$$
Introducing the $(m-1)$-forms
$$\omega_\la =\di \frac{\pa}{\pa t^\la}\rfloor\omega,$$
a {\it costate variable vector or Lagrange multiplier vector} $p=p^\al(t)\di\frac{\pa}{\pa t^\al}$ is identified to
the $(m-1)$-form $p = p^\la(t) \omega_\la$. We build a {\it Hamiltonian $m$-form}
$${\cal H} = {\cal L}(t,x(t),u(t)) + p^\la(t) a_i(x(t), u(t))dx^i(t)\wedge \omega_\la.$$

{\bf Theorem} ({\bf Multitime maximum principle})
{\it Suppose that the problem of maximizing the functional (11) constrained by (12) has
an interior optimal solution $\hat u(t)$, which determines the optimal evolution $x(t)$. Then there exists a
costate function $(p(t))$ such that
$$
\frac{\partial {\mathcal H}}{\partial p^\la}= a_i(x(t), u(t))dx^i(t)\wedge \omega_\la = 0, \leqno(13)
$$
the function $(p(t))$ is the unique solution of the following Pfaff system (adjoint system)
$$
d(p^\la a_j(x,u))\wedge \omega_\la = \frac{\pa {\cal H}}{\pa x^j} \leqno(14)
$$
and satisfies the critical point conditions
$$
{\mathcal H}_{u^a}\left( t,x(t),u(t),p(t)\right) =0,\,a = \overline{1,k}.\leqno(15)
$$}

{\bf Proof} We use the Hamiltonian $m$-form ${\cal H}$.
The solutions of the foregoing problem are among the solutions of the free maximization problem of the functional
$$J(u(\cdot)) = \int_{\tilde\Omega}{\cal H}(t,x(t),u(t),p(t)) + \int_{\pa {\tilde\Omega}} g(x(t))d\si\,$$
where $\tilde \Omega = (\Omega_{t_0t_1}, x(\Omega_{t_0t_1}))\subset R^m_+\times R^n$.

Suppose that there exists a continuous control $\hat u(t)$ defined over the interval
$\Omega_{t_0t_1}$ with $\hat u(t)\in \hbox{Int}\,{U}$
which is an optimum point in the previous problem. Now we consider a control variation
$u(t,\epsilon)=\hat u(t) + \epsilon h(t)$, where $h$ is an arbitrary continuous vector function,
and a state variation $x(t,\ep),\, t \in \Omega_{t_0t_1}$, connected by
$$a_i(x(t;\epsilon), u(t; \epsilon))dx^i(t;\epsilon) = 0,\, x(t;0) = x(t),\, x(t_0,0)=x_0.$$
Since $\hat u(t)\in \hbox{Int}\,{\cal U}$ and any continuous function over a compact set
$\Omega_{t_0t_1}$ is bounded, there exists $\epsilon_h>0$ such
that $u(t,\epsilon)=\hat u(t) + \epsilon h(t)\in \hbox{Int}\,{U},\,\,\forall |\epsilon|<\epsilon_h$.
This $\epsilon$ is used in our variational arguments.

For $|\epsilon|<\epsilon_h$, we define the function
$$J(\epsilon)= \int_{\tilde\Omega(\ep)}{\cal H}(t,x(t,\epsilon),u(t,\epsilon),p(t)) + \int_{\pa {\tilde\Omega(\ep)}}g(x(t,\ep))d\si$$
$$= \int_{\tilde\Omega(\ep)}{\cal L}(t,x(t,\epsilon),u(t,\epsilon)) + p^\la(t) a_i(x(t,\epsilon), u(t,\epsilon)) dx^i(t,\epsilon)\wedge \omega_\la +
\int_{\pa {\tilde\Omega(\ep)}}g(x(t,\ep))d\si.$$
Differentiating with respect to $\epsilon$, it follows
$$J^{\prime}(\epsilon)= \int_{\tilde\Omega(\ep)}{\cal L}_{x^j}(t,x(t,\epsilon),u(t,\epsilon))x^j_{\epsilon}(t,\epsilon) +
\int_{\pa {\tilde\Omega(\ep)}}g_{x^j}(x(t,\ep))x^j_{\epsilon}(t,\epsilon) d\si$$
$$ + \int_{\tilde\Omega(\ep)} p^\la(t) \frac{\pa a_i}{\pa x^j}(x(t,\epsilon), u(t,\epsilon))x^j_{\epsilon}(t,\epsilon) dx^i(t,\epsilon)\wedge \omega_\la$$
$$ +  \int_{\tilde\Omega(\ep)} p^\la(t) a_i(x(t,\epsilon), u(t,\epsilon)) dx^i_{\epsilon}(t,\epsilon)\wedge \omega_\la$$
$$ + \int_{\tilde\Omega(\ep)}{\cal L}_{u^a}(t,x(t,\epsilon),u(t,\epsilon)h^a(t)$$
$$+ \int_{\tilde\Omega(\ep)} p^\la(t) \frac{\pa a_i}{\pa u^a}(x(t,\epsilon), u(t,\epsilon))h^a(t) dx^i(t,\epsilon))\wedge \omega_\la.$$
Evaluating at $\epsilon = 0$, we find $\tilde\Omega(0) = \tilde\Omega$ and
$$J^{\prime}(0)= \int_{\tilde\Omega}\left({\cal L}_{x^j}(t,x(t),\hat u(t))+ p^\la(t) \frac{\pa a_i}{\pa x^j}(x(t), \hat u(t))dx^i(t)\wedge \omega_\la\right)x^j_{\epsilon}(t,0)$$
$$ +  \int_{\tilde\Omega} p^\la(t) a_i(x(t), \hat u(t)) dx^i_{\epsilon}(t,0)\wedge \omega_\la +\int_{\pa {\tilde\Omega}}g_{x^j}(x(t))x^j_{\epsilon}(t,0) d\si$$
$$ + \int_{\tilde\Omega}\left({\cal L}_{u^a}(t,x(t),\hat u(t)) + p^\la(t) \frac{\pa a_i}{\pa u^a}(x(t), \hat u(t))dx^i(t)\wedge \omega_\la\right) h^a(t),$$
where $x(t)$ is the $m$-sheet of the state variable corresponding to the optimal control $\hat u(t)$.

To evaluate the multiple integral

$$\int_{\tilde\Omega} p^\la(t) a_i(x(t), \hat u(t)) dx^i_{\epsilon}(t,0)\wedge \omega_\la,$$
we integrate by parts, via the formula
$$d(p^\la a_i x^i_{\epsilon} \omega_\la)=(x^i_\epsilon d(p^\la a_i) + p^\la a_i dx^i_\epsilon)\wedge\omega_\la,$$
obtaining
$$\di\int_{\tilde\Omega}p^\la a_idx^i_\epsilon\wedge \omega_\la=\di\int_{\tilde\Omega}d(p^\la a_ix^i_\epsilon\omega_\la)
- \di\int_{\tilde\Omega}d(p^\la a_i) x^i_\epsilon\wedge \omega_\la,$$
where the symbol $d$ in the last integral means the differentiation with respect to $p$ and $x$.
Now we apply the Stokes integral formula
$$\di\int_{\tilde\Omega}d(p^\la a_ix^i_\epsilon\omega_\la) = \di\int_{\pa{\tilde\Omega}}\delta_{\al\be}p^\al a_ix^i_\epsilon n^\be d\sigma,$$
where $(n^\be(t))$ is the unit normal vector to the boundary ${\pa{\tilde\Omega}}$.
Since the integral from the middle can be written
$$\di\int_{\tilde\Omega}p^\la a_idx^i_\epsilon\wedge \omega_\la=\di\int_{\pa{\tilde\Omega}}\delta_{\al\be}p^\al a_ix^i_\epsilon n^\be3 d\sigma
- \di\int_{\tilde\Omega}d(p^\la a_i) x^i_\epsilon\wedge \omega_\la,$$
we find $J^{\prime}(0)$ as
$$\int_{\tilde\Omega}{\cal L}_{x^j}(t,x(t),\hat u(t))x^j_{\epsilon}(t,0)$$
$$+ \int_{\tilde\Omega} \left(p^\la(t) \frac{\pa a_i}{\pa x^j}(x(t), \hat u(t))dx^i(t)
- d(p^\la a_j(x(t), \hat u(t))\right)\wedge \omega_\la x^j_{\epsilon}(t,0)$$
$$+ \di\int_{\pa{\tilde\Omega}}\delta_{\al\be}p^\al a_ix^i_\epsilon(t,0) n^\be(t) d\sigma +\int_{\pa {\tilde\Omega}}g_{x^j}(x(t))x^j_{\epsilon}(t,0) d\si$$
$$ + \int_{\tilde\Omega}\left({\cal L}_{u^a}(t,x(t),\hat u(t)) + p^\la(t) \frac{\pa a_i}{\pa u^a}(x(t), \hat u(t))dx^i(t)\wedge \omega_\la\right) h^a(t).$$
We select the costate function $p(t)$ as solution of the adjoint Pfaff equation (boundary value problem)
$${\cal L}_{x^j}(t,x(t),\hat u(t)) + \left(p^\la(t) \frac{\pa a_i}{\pa x^j}(x(t), \hat u(t))dx^i(t)
- d(p^\la a_j(x(t), \hat u(t))\right)\wedge \omega_\la = 0,$$
$$\left(\delta_{\al\be}p^\al(t) n^\be(t) a_i(x(t), \hat u(t)) + g_{x^i}(x(t))\right)|_{\pa \Omega} = 0.$$
On the other hand, we need $J^{\prime}(0)=0$ for all $h(t)=(h^a(t))$. Since
the variation $h$ is arbitrary, we get (critical point condition)
$$\frac{\pa {\cal L}}{\pa {u^a}}(t,x(t),\hat u(t)) + p^\la(t) \frac{\pa a_i}{\pa u^a}(x(t), \hat u(t))dx^i(t)\wedge \omega_\la = 0.$$

The foregoing equations (14) and (15) can be written
$$\frac{\pa {\cal H}}{\pa {x^j}} - d(p^\la a_j) \wedge \omega_\la = 0,\,\,\frac{\pa {\cal H}}{\pa {u^a}} = 0.$$

Let us start with a

\subsubsection{Multitime optimal control problem of \\maximizing a curvilinear integral functional}

{\it Find}
{$$\max_{u(\cdot)}\,I(u(\cdot)) = \int_{\Gamma_{t_0t_1}}L_\al(t,x(t),u(t)) dt^\al + g(x(t_1))\leqno(16)$$
{\it subject to}
$$a_i(x(t), u(t))dx^i(t) = 0,\,\hbox{a.e.}\,\,\, t\in \Omega_{t_0t_1},\, x(t_0) = x_0.\leqno(17)$$}
It is supposed that $L_\al: \Omega_{t_0t_1} \times A\times U \to R$ and $a_i: A\times U \to R, \, i = 1,..., n$
are $C^2$ functions. Ingredients: ${\cal L} = L_\al(t,x(t),u(t)) dt^\al$ is an $1$-form, $A$ is a bounded and closed subset of $R^n$,
containing the images of the $m$-sheets $x(t), \,t\in \Omega_{t_0t_1}$  of the controlled system,
and $x_0$ and $x_1$ are the initial and final states of the $m$-sheet $x(t)$ in the controlled system.
The set, in which the control functions $u^a$ takes their values, is called as $U$, which is a bounded and
closed subset of $R^k$.

Let us find the first order necessary conditions for an optimal pair $(x,u)$.
We fix the control $u(t)$ and variate the state $x(t)$ into $x(t,\ep),\, \ep = (\ep^1,...,\ep^m)$.
We find the {\it multitime infinitesimal deformation (Pfaff) system}
$$\frac{\pa a_i}{\pa x^j}(x,u)y^j_\al(t) dx^i + a_i(x,u) dy^i_\al = 0$$
of the nonholonomic constraint $a_i(x(t), u(t))dx^i(t) = 0$.
It follows the {\it multitime adjoint Pfaff system}
$$d(pa_j(x,u)) = \frac{\pa (pa_idx^i)}{\pa x^j},$$
whose solution $p(t)$ is called the {\it costate vector}. Here, the symbol $d$ in the left hand member of the adjoint equation
means the differentiation with respect to $p$ and $x$.

We use the Lagrangian $1$-form
$${\cal L} = L_\al(t,x(t),u(t)) dt^\al.$$
Introducing a {\it costate variable or Lagrange multiplier} $p$, we build a {\it Hamiltonian $1$-form}
$${\cal H} = {\cal L}(t,x(t),u(t)) + p(t) a_i(x(t), u(t))dx^i(t).$$

{\bf Theorem} ({\bf Multitime maximum principle})
{\it Suppose that the problem of maximizing the functional (16) constrained by (17) has
an interior optimal solution $\hat u(t)$, which determines the optimal evolution $x(t)$. Then there exists a
costate function $(p(t))$ such that
$$
\frac{\partial {\mathcal H}}{\partial p}= a_i(x(t), u(t))dx^i(t) = 0, \leqno(18)
$$
the function $(p(t))$ is the unique solution of the following Pfaff system (adjoint system)
$$
d(pa_j(x,u)) = \frac{\pa {\cal H}}{\pa x^j}\leqno(19)
$$
and the critical point conditions
$$
{\mathcal H}_{u^a}\left( t,x(t),u(t),p(t)\right) =0,\,a = \overline{1,k}\leqno(20)
$$
hold.}

{\bf Proof} We use the Hamiltonian $1$-form ${\cal H}$.
The solutions of the foregoing problem are between the solutions of the free maximization problem of the curvilinear integral functional
$$J(u(\cdot)) = \int_{\tilde\Gamma}{\cal H}(t,x(t),u(t)) + g(x(t_1)),$$
where $\tilde \Gamma = (\Gamma_{t_0t_1}, x(\Gamma_{t_0t_1}))\subset R^m_+\times R^n$.


Suppose that there exists a continuous control $\hat u(t)$ defined over the interval $\Omega_{t_0t_1}$ with $\hat u(t)\in \hbox{Int}\,{U}$
which is an optimum point in the previous problem. We consider a control variation
$u(t,\epsilon)=\hat u(t) + \epsilon h(t)$, where $h$ is an arbitrary continuous vector function,
and a state variation $x(t,\ep),\, t \in \Omega_{t_0t_1}$, connected by
$$a_i(x(t;\epsilon), u(t; \epsilon))dx^i(t;\epsilon) = 0,\, x(t;0) = x(t),\, x(t_0,0)=x_0.$$
Since $\hat u(t)\in \hbox{Int}\,{\cal U}$ and a continuous function over a compact set
$\Omega_{t_0t_1}$ is bounded, there exists $\epsilon_h>0$ such
that $u(t,\epsilon)=\hat u(t) + \epsilon h(t)\in \hbox{Int}\,{U},\,\,\forall |\epsilon|<\epsilon_h$.
This $\epsilon$ is used in our variational arguments.

For $|\epsilon|<\epsilon_h$, we define the function
$$J(\epsilon)= \int_{\tilde\Gamma(\ep)}{\cal H}(t,x(t,\epsilon),u(t,\epsilon),p(t)) + g(x(t_1,\ep))$$
$$= \int_{\tilde\Gamma(\ep)}{\cal L}(t,x(t,\epsilon),u(t,\epsilon)) + p(t) a_i(x(t,\epsilon), u(t,\epsilon)) dx^i(t,\epsilon)+ g(x(t_1,\ep)).$$
Differentiating with respect to $\epsilon$, it follows
$$
J^{\prime}(\epsilon)= \int_{\tilde\Gamma(\ep)}{\cal L}_{x^j}(t,x(t,\epsilon),u(t,\epsilon))x^j_{\epsilon}(t,\epsilon)+ g_{x^j}(x(t_1,\ep))x^j_{\epsilon}(t,\epsilon)
$$
$$ + \int_{\tilde\Gamma(\ep)} p(t) \frac{\pa a_i}{\pa x^j}(x(t,\epsilon), u(t,\epsilon))x^j_{\epsilon}(t,\epsilon) dx^i(t,\epsilon)$$
$$ +  \int_{\tilde\Gamma(\ep)} p(t) a_i(x(t,\epsilon), u(t,\epsilon)) dx^i_{\epsilon}(t,\epsilon)$$
$$ + \int_{\tilde\Gamma(\ep)}{\cal L}_{u^a}(t,x(t,\epsilon),u(t,\epsilon))h^a(t)$$
$$+ \int_{\tilde\Gamma(\ep)} p(t) \frac{\pa a_i}{\pa u^a}(x(t,\epsilon), u(t,\epsilon))h^a(t) dx^i(t,\epsilon)).$$
Evaluating at $\epsilon = 0$, we find $\tilde\Gamma(0) = \tilde\Gamma$ and
$$J^{\prime}(0)= \int_{\tilde\Gamma}\left({\cal L}_{x^j}(t,x(t),\hat u(t))+ p(t) \frac{\pa a_i}{\pa x^j}(x(t), \hat u(t))dx^i(t)\right)x^j_{\epsilon}(t,0)$$
$$
+  \int_{\tilde\Gamma} p(t) a_i(x(t), \hat u(t)) dx^i_{\epsilon}(t,0) + g_{x^j}(x(t_1)x^j_{\epsilon}(t,0)
$$
$$ + \int_{\tilde\Gamma}\left({\cal L}_{u^a}(t,x(t),\hat u(t)) + p(t) \frac{\pa a_i}{\pa u^a}(x(t), \hat u(t))dx^i(t)\right) h^a(t),$$
where $x(t)$ is the $m$-sheet of the state variable corresponding to the optimal control $\hat u(t)$.

To evaluate the curvilinear integral
$$\int_{\tilde\Gamma} p(t) a_i(x(t), \hat u(t)) dx^i_{\epsilon}(t,0),$$
we integrate by parts, via the formula
$$d(p a_i x^i_{\epsilon})= x^i_\epsilon d(p a_i) + p a_i dx^i_\epsilon,$$
obtaining
$$\di\int_{\tilde\Gamma}p a_idx^i_\epsilon = \di\int_{\tilde\Gamma}d(p a_ix^i_\epsilon)
- \di\int_{\tilde\Gamma}d(p a_i) x^i_\epsilon$$
$$ = (p(t) a_i(x(t), \hat u(t)) x^i_{\epsilon}(t,0))|_{t_0}^{t_1} - \di\int_{\tilde\Gamma}d(p(t) a_i(x(t), \hat u(t))) x^i_\epsilon,$$
where the symbol $d$ in the last integral means the differentiation with respect to $p$ and $x$. We find
$$J^{\prime}(0) = \int_{\tilde\Gamma}{\cal L}_{x^j}(t,x(t),\hat u(t))x^j_{\epsilon}(t,0)$$
$$+ \int_{\tilde\Gamma} \left(p(t)\frac{\pa a_i}{\pa x^j}(x(t), \hat u(t))dx^i(t) - d(p(t)a_j(x(t), \hat u(t))\right) x^j_{\epsilon}(t,0)$$
$$+ (p(t) a_i(x(t),\hat u(t))x^i_{\epsilon}(t,0))|_{t_0}^{t_1} + g_{x^j}(x(t_1)x^j_{\epsilon}(t,0)$$
$$ + \int_{\tilde\Gamma}\left({\cal L}_{u^a}(t,x(t),\hat u(t)) + p(t) \frac{\pa a_i}{\pa u^a}(x(t), \hat u(t))dx^i(t)\right) h^a(t).$$
We select the costate function $p(t)$ as solution of the adjoint Pfaff equation (terminal value problem)
$${\cal L}_{x^j}(t,x(t),\hat u(t)) + p(t)\frac{\pa a_i}{\pa x^j}(x(t), \hat u(t))dx^i(t) - d(p(t)a_j(x(t), \hat u(t)) = 0,$$
subject to $p(t_1) a_i(x(t_1),\hat u(t_1)) + g_{x^i}(x(t_1)) =0$.
On the other hand, we need $J^{\prime}(0)=0$ for all $h(t)=(h^a(t))$. Since
the variation $h$ is arbitrary, we get (critical point condition)
$$\frac{\pa {\cal L}}{\pa {u^a}}(t,x(t),\hat u(t)) + p(t) \frac{\pa a_i}{\pa u^a}(x(t), \hat u(t))dx^i(t) = 0.$$

The foregoing equations (19) and (20) can be written
$$\frac{\pa {\cal H}}{\pa {x^j}} = d(p\, a_j),\,\,\frac{\pa {\cal H}}{\pa {u^a}} = 0.$$

\section{Optimal control on distributions \\described by vector fields}

\subsection{Infinitesimal deformations and adjointness on \\distributions}

The same distribution $D$ can be described in terms of smooth vector fields (or generators),
$$D = \hbox{span}\{X_a(x)|\, a_i(x)X^i_a = 0,\, a = 1,...,n-1\},\leqno (16)$$
if and only if $n\geq 3$. Any vector field $Y$ in $D$ can be written in the form $Y(x) = u^a(x)X_a(x)$.

Let $x(t)$ be a curve solution of the differential system
$$\dot x(t) = u^a(x(t))X_a(x(t)).$$
Let $x(t;\epsilon)$ be a differentiable variation of $x(t)$, i.e.,
$$\dot x(t;\epsilon) = u^a(x(t;\epsilon))X_a(x(t;\epsilon)),\, x(t;0) = x(t).$$ Denoting $y^i(t) = \di\frac{\pa x^i}{\pa \epsilon}(t;0)$,
we find the {\it single-time infinitesimal deformation system}
$$\dot y^j(t) = \left(\frac{\pa u^a}{\pa x^i}(x(t)) X^j_a(x(t)) + u^a(x(t)) \frac{\pa X^j_a}{\pa x^i}(x(t))\right)y^i(t).\leqno(17)$$
The {\it single-time adjoint (dual) system} is
$$\dot p_k (t)= - \left(\frac{\pa u^a}{\pa x^k}(x(t)) X^j_a(x(t)) + u^a(x(t)) \frac{\pa X^j_a}{\pa x^k}(x(t))\right)p_j(t),\leqno(18)$$
whose solution $p=(p_k)$ is called the {\it costate vector}.
The foregoing PDE systems (17) and (18) are {\it adjoint} ({\it dual}) in the sense of
{\it constant interior product of solutions}, i.e., the scalar product $p_k\,y^k$ is a first integral.

Let $x(t)$ be an $m$-sheet integral submanifold of the distribution $D$, i.e., a solution of the multitime
partial differential system $$\frac{\pa x}{\pa t^\al}(t) = u^a_\al(x(t))X_a(x(t)),\, \al = 1,...,m < n-1.$$
Let $\epsilon = (\epsilon^\al),\,\al = 1,...,m$ and let $x(t;\epsilon)$ be a differentiable variation of $x(t)$, i.e.,
$$\frac{\pa x}{\pa t^\al}(t;\epsilon) = u^a_\al(x(t;\epsilon))X_a(x(t;\epsilon)),\,x(t;0) = x(t).$$
Introducing the vector fields $y^i_\al(t) = \di\frac{\pa x^i}{\pa \epsilon^\al}(t;0)$,
we find the {\it multitime infinitesimal deformation system}
$$\frac{\pa y^j_\al}{\pa t^\be}(t) = \left(\frac{\pa u^a_\be}{\pa x^i}(x(t)) X^j_a(x(t)) + u^a_\be(x(t)) \frac{\pa X^j_a}{\pa x^i}(x(t))\right)y^i_\al(t).\leqno(19)$$
The {\it multitime adjoint (dual) system} is
$$\frac{\pa p_k^\al}{\pa t^\be} (t)= - \left(\frac{\pa u^a_\be}{\pa x^k}(x(t)) X^j_a(x(t)) + u^a_\be(x(t)) \frac{\pa X^j_a}{\pa x^k}(x(t))\right)p_j^\al(t),\leqno(20)$$
whose solution $p=(p^\al_k)$ is called the {\it costate matrix}.
The foregoing PDE systems (19) and (20) are {\it adjoint} ({\it dual}) in the sense of
{\it constant interior product of solutions}, i.e., the scalar product $p^\al_k\,y^k_\al$ is a first integral.

Of course, taking the trace, we can define the {\it costate matrix} $p: \Omega_{0T} \to R^{mn},\,\,p=(p^\al_k),$
as the solution of the {\it divergence adjoint PDE system} (trace)
$$\frac{\pa p_k^\al}{\pa t^\al} (t)= - \left(\frac{\pa u^a_\al}{\pa x^k}(x(t)) X^j_a(x(t)) + u^a_\al(x(t)) \frac{\pa X^j_a}{\pa x^k}(x(t))\right)p_j^\al(t).\leqno(21)$$
But than, the PDEs systems (19) and (20) are {\it adjoint (dual)} in the sense of {\it zero total divergence}
of the tensor field $Q^\al_\be  = p^\al_k\,y^k_\be$ produced by their solutions.
The divergence dual PDE system (21) has solutions since it contains $n$ PDEs with $nm$ unknown functions $p^\al_i $.
We can select a solution of the gradient form $p^\al_k(t) = \di\frac{\pa v^\al}{\pa x^k}(t,x(t))$.

{\bf Remark} The {\it multitime adjoint Pfaff system} can be defined independent on the dimension
of the parameter $\ep$. Particularly, the {\it multitime adjoint Pfaff system} can be
$$\frac{\pa p_k}{\pa t^\be} (t)= - \left(\frac{\pa u^a_\be}{\pa x^k}(x(t)) X^j_a(x(t)) + u^a_\be(x(t)) \frac{\pa X^j_a}{\pa x^k}(x(t))\right)p_j(t).\leqno(*)$$

\subsection{Single-time optimal control problems on \\a distribution}
Let
$$D = \hbox{span}\{X_a(x)|\, a_i(x)X^i_a = 0,\, a = 1,...,n-1\}.$$
be a distribution on $R^n$ and $x(t),\,\,t\in I = [t_0, t_1],$ be an integral curve of the driftless control system
$$dx(t) = u^a(x(t))X_a(x(t))dt.$$

A {\it single-time optimal control problem is defined to be maximizing the functional
$$I(u(\cdot)) = \int_{t_0}^{t_1}L(t,x(t),u(x(t))) dt\leqno(22)$$
subject to
$$dx(t) = u^a(x(t))X_a(x(t))dt,\,\hbox{a.e.}\,\,\, t\in I = [t_0, t_1],\, x(t_0) = x_1,\, x(t_1)= x_1.\leqno(23)$$}
It is supposed that $L: I \times A\times U \to R$ is a $C^2$ function and $X_a: A \to R, \,a = 1,..., n-1$
are $C^2$ functions. Ingredients: $A$ is a bounded and closed subset of $R^n$, which the trajectory of controlled system is constrained to stay
for $t\in I$, and $x_0$ and $x_1$ are the initial and final states of the trajectory $x(t)$ in the controlled system.
The set in which the control functions $u^a$ takes their values in it, is called as $U$, which is a bounded and
closed subset of $R^{n-1}$. The map $u$ is assumed to be piecewise smooth or piecewise analytic. Such maps are
called {\it admissible} and the space ${\cal U}$ of all such maps is called the {\it set of admissible controls}.

Let us find the first order necessary conditions for an optimal pair $(x,u)$. Firstly, the {\it single-time infinitesimal deformation (Pfaff) equation}
of the constraint $dx(t) = u^a(x(t))X_a(x(t))dt$ is the system (17) the {\it single-time adjoint Pfaff equation} is the system (18).

The control variables may be {\it open-loop} $u^a(t)$, depending directly on the time variable $t$, or {\it closed-loop}
(or {\it feedback}) $u^a(x(t))$, depending on the state $x(t)$.

{\bf Open-loop control variables}

To simplify, we accept an open-loop control $u^a(t)$. Using the Lagrangian $1$-form
$${\cal L}(t,x(t),u(t),p(t)) = L(t,x(t),u(t))dt+p_i(t)[u^a(t)X^i_a(x(t))dt - dx^i(t)],$$
we build the {\it Hamiltonian $1$-form}
$${\cal H} = L(t,x(t),u(t)) dt+ p_i(t) u^a(t) X^i_a(x(t))dt.$$

{\bf Theorem} ({\bf Single-time maximum principle})
{\it Suppose that the problem of maximizing the functional (22) constrained by (23) has
an interior optimal solution $\hat u(t)$, which determines
the optimal evolution $x(t)$. Then there exists a
costate vector p(t) = $(p_i(t))$ such that
$$
dx^i = \frac{\partial {\mathcal H}}{\partial p_i},\leqno(24)
$$
the function $p(t)$ is the unique solution of the following Pfaff system (adjoint system)
$$
dp_i = - \frac{\pa {\cal H}}{\pa x^i}\leqno(25)
$$
and the critical point conditions
$$
{\mathcal H}_{u^a}\left( t,x(t),u(t),p(t)\right) =0,\,a = \overline{1,n-1}\leqno(26)
$$
hold.}

{\bf Proof} We use the Lagrangian $1$-form ${\cal L}$.
The solutions of the forgoing problem are between the solutions of the free maximization problem of the curvilinear integral functional
$$J(u(\cdot)) = \int_{\tilde\Gamma}{\cal L}(t,x(t),u(t),p(t)),$$
where $\tilde \Gamma = ([t_0,t_1], x([t_0,t_1]))\subset R_+\times R^n$.

Suppose that there exists a continuous control $\hat u(t)$ defined over the interval
$I$ with $\hat u(t)\in \hbox{Int}\,{U}$
which is an optimum point in the previous problem. Now consider a variation
$u(t,\epsilon)=\hat u(t) + \epsilon h(t)$, where $h$ is an arbitrary continuous vector function.
Since $\hat u(t)\in \hbox{Int}\,{\cal U}$ and a continuous function over a compact set
$I$ is bounded, there exists $\epsilon_h>0$ such
that $u(t,\epsilon)=\hat u(t) + \epsilon h(t)\in \hbox{Int}\,{U},\,\,\forall |\epsilon|<\epsilon_h$.
This $\epsilon$ is used in our variational arguments.

Define $x(t,\epsilon)$ as the $1$-sheet of the state variable corresponding to the control variable $u(t,\epsilon)$, i.e.,
$$dx^i(t;\epsilon) = u^a(x(t;\epsilon))X^i_a(x(t;\epsilon))dt,\, x(t;0) = x(t).$$
and $x(t_0,0)=x_0$, $x(t_1,0)=x_1$. For $|\epsilon|<\epsilon_h$, we define the function
$$J(\epsilon)= \int_{\tilde\Gamma(\ep)}{\cal L}(t,x(t,\epsilon),u(t,\epsilon),p(t))$$
$$ \int_{\tilde\Gamma(\ep)} L(t,x(t,\epsilon),u(t,\epsilon))dt+p_i(t)[u^a(t,\epsilon)X^i_a(x(t,\epsilon))dt - dx^i(t,\epsilon)].$$
Differentiating with respect to $\epsilon$, it follows
$$J^{\prime}(\epsilon)= \int_{\tilde\Gamma(\ep)}\left(L_{x^j}(t,x(t,\epsilon),u(t,\epsilon)) + p_i(t)u^a(t,\epsilon) X^i_{ax^j}\right)x^j_{\epsilon}(t,\epsilon)dt$$
$$ -  \int_{\tilde\Gamma(\ep)} p_i(t) dx^i_{\epsilon}(t,\epsilon)$$
$$ + \int_{\tilde\Gamma(\ep)}\left({L}_{u^a}(t,x(t,\epsilon),u(t,\epsilon) + p_i(t) X^i_a(x(t,\epsilon))\right)h^a(t) dt.$$
Evaluating at $\epsilon = 0$, we find $\tilde\Gamma(0) = \tilde\Gamma$ and
$$J^{\prime}(0)= \int_{\tilde\Gamma}\left(L_{x^j}(t,x(t),\hat u(t)) + p_i(t)\hat{u}^a(t) X^i_{ax^j}\right)x^j_{\epsilon}(t,0)dt$$
$$ -  \int_{\tilde\Gamma} p_i(t) dx^i_{\epsilon}(t,0)$$
$$ + \int_{\tilde\Gamma}\left({L}_{u^a}(t,x(t),\hat u(t) + p_i(t) X^i_a(x(t))\right)h^a(t) dt,$$
where $x(t)$ is the curve of the state variable corresponding to the optimal control $\hat u(t)$. Since the integral from the middle can be written
$$\int_{\tilde\Gamma} p_i(t) dx^i_{\epsilon}(t,0) = p_i(t) x^i_{\epsilon}|_{\pa\tilde\Gamma} - \int_{\tilde\Gamma}x^i_{\epsilon}(t,0)dp_i(t),$$
we find $J^{\prime}(0)$ as
$$J^{\prime}(0)= \int_{\tilde\Gamma}\left(\left(L_{x^j}(t,x(t),\hat u(t)) + p_i(t)\hat{u}^a(t) X^i_{ax^j}\right)dt + dp_j\right)x^j_{\epsilon}(t,0)$$
$$ -  p_i(t) x^i_{\epsilon}(t,0)|_{t_0}^{t_1}$$
$$ + \int_{\tilde\Gamma}\left({L}_{u^a}(t,x(t),\hat u(t) + p_i(t) X^i_a(x(t))\right)h^a(t) dt,$$
We select the costate function $p(t)$ as solution of the adjoint Pfaff equation (boundary value problem)
$$\left(L_{x^j}(t,x(t),\hat u(t)) + p_i(t)\hat{u}^a(t) X^i_{ax^j}\right)dt + dp_j = 0,\,\,p(t_0) = 0, p(t_1) = 0.$$
On the other hand, we need $J^{\prime}(0)=0$ for all $h(t)=(h^a(t))$. Since
the variation $h$ is arbitrary, we get (critical point condition)
$${L}_{u^a}(t,x(t),\hat u(t) + p_i(t) X^i_a(x(t)) = 0.$$

{\bf Example} Consider the ODE system $\dot x^1(t) = {x^2}^2(t), \dot x^2(t) = u(t)$
generated by the vector fields $X = {x^2}^2 \frac{\pa}{\pa x^1}$, $Y = \frac{\pa}{\pa x^2}$.
We compute the Lie brackets
$$[X,Y] = - 2 x^2 \frac{\pa}{\pa x^1}, [[X,Y],Y] = 2 \frac{\pa}{\pa x^1}.$$
The vector fields $Y$ and $[[X,Y],Y]$ are linearly independent.
On the other hand, the $x^1$-coordinate is increasing since ${x^2}^2 \geq 0$.
Consequently, the system is not really controllable.


{\bf Closed-loop control variables}

Now, we accept a closed-loop control $u^a(x(t))$. Using the Lagrangian $1$-form
$${\cal L} = L(t,x(t),u(x(t)))dt+p_i(t)[u^a(x(t))X^i_a(x(t))dt - dx^i(t)],$$
we build the {\it Hamiltonian $1$-form}
$${\cal H} = L(t,x(t),u(x(t))dt + p_i(t) u^a(x(t)) X^i_a(x(t))dt.$$

{\bf Theorem} ({\bf Single-time maximum principle})
{\it Suppose that the problem of maximizing the functional (22) constrained by (23) has
an interior optimal solution $\hat u(x(t))$, which determines
the optimal evolution $x(t)$. Then there exists a
costate vector p(t) = $(p_i(t))$ such that
$$
dx^i = \frac{\partial {\mathcal H}}{\partial p_i},
$$
the function $p(t)$ is the unique solution of the following Pfaff system (adjoint system)
$$
dp_i = - \frac{\pa {\cal H}}{\pa x^i}
$$
and the critical point conditions
$$
{\mathcal H}_{u^a}\left( t,x(t),u(x(t)),p(t)\right) =0,\,a = \overline{1,n-1}
$$
hold.}

{\bf Proof} The new functional is
$$J(x(\cdot), u(x(\cdot))) = \int_{\tilde\Gamma}{\cal L}(t,x(t),u(x(t)),p(t)).$$
A variation $x(t,\ep)$ induces a variation $u(x(t,\ep)) =\hat u(x(t)) + \ep h(t)$. Then
$$J(\ep) = \int_{\tilde\Gamma(\ep)}{\cal L}(t,x(t,\ep), u(x(t,\ep)), p(t)).$$
It follows
$$
J^{\prime}(\epsilon)= \int_{\tilde\Gamma(\ep)}L_{x^j}(t,x(t,\epsilon),u(x(t,\epsilon)) x^j_{\epsilon}(t,\epsilon)dt
$$
$$
+ \int_{\tilde\Gamma(\ep)}\left(p_i(t)u^a_{x^j}X^i_a+ p_i(t)u^a(x(t,\epsilon)) X^i_{ax^j}\right)x^j_{\epsilon}(t,\epsilon)dt
$$
$$ -  \int_{\tilde\Gamma(\ep)} p_i(t) dx^i_{\epsilon}(t,\epsilon)$$
$$ + \int_{\tilde\Gamma(\ep)}\left({L}_{u^a}(t,x(t,\epsilon),u(x(t,\epsilon)) + p_i(t) X^i_a(x(t,\epsilon))\right)h^a(t) dt.$$
Evaluating at $\epsilon = 0$, we find $\tilde\Gamma(0) = \tilde\Gamma$ and
$$
J^{\prime}(0)= \int_{\tilde\Gamma}L_{x^j}(t,x(t),\hat u(x(t)) x^j_{\epsilon}(t,0)dt
$$
$$
+ \int_{\tilde\Gamma}\left(p_i(t)\hat u^a_{x^j}X^i_a+ p_i(t)u^a(x(t)) X^i_{ax^j}\right)x^j_{\epsilon}(t,0)dt
$$
$$ -  \int_{\tilde\Gamma} p_i(t) dx^i_{\epsilon}(t,0)$$
$$ + \int_{\tilde\Gamma}\left({L}_{u^a}(t,x(t),\hat u(x(t))) + p_i(t) X^i_a(x(t))\right)h^a(t) dt.$$
Since the integral from the middle can be written
$$\int_{\tilde\Gamma} p_i(t) dx^i_{\epsilon}(t,0) = p_i(t) x^i_{\epsilon}|_{\pa\tilde\Gamma} - \int_{\tilde\Gamma}x^i_{\epsilon}(t,0)dp_i(t),$$
we find $J^{\prime}(0)$ as
$$
J^{\prime}(0)= \int_{\tilde\Gamma}\left(\left(L_{x^j}(t,x(t),\hat u(x(t))) + p_i(t)\hat u^a_{x^j}X^i_a +p_i(t)\hat{u}^a(x(t)) X^i_{ax^j}\right)dt + dp_j\right)x^j_{\epsilon}(t,0)
$$
$$ -  p_i(t) x^i_{\epsilon}(t,0)|_{t_0}^{t_1}$$
$$ + \int_{\tilde\Gamma}\left({L}_{u^a}(t,x(t),\hat u(x(t)) + p_i(t) X^i_a(x(t))\right)h^a(t) dt,$$
We select the costate function $p(t)$ as solution of the adjoint Pfaff equation (boundary value problem)
$$\left(L_{x^j}(t,x(t),\hat u(x(t))) + p_i(t)\hat{u}^a(x(t)) X^i_{ax^j}\right)dt + dp_j = 0,\,\,p(t_0) = 0, p(t_1) = 0.$$
On the other hand, we need $J^{\prime}(0)=0$ for all $h(t)=(h^a(t))$. Since
the variation $h$ is arbitrary, we get (critical point condition)
$${L}_{u^a}(t,x(t),\hat u(x(t)) + p_i(t) X^i_a(x(t)) = 0.$$

\subsection{Multitime optimal control problems on \\a distribution}

Let
$$D = \hbox{span}\{X_a(x)|\, a_i(x)X^i_a = 0,\, a = 1,...,n-1\},$$
$n\geq 3$, be a distribution on $R^n$ and $x(t),\,\,t\in \Omega_{t_0t_1}$, be an $m$-sheet of the driftless control system
$$dx(t) = u^a_\al(x(t))X_a(x(t))dt^\al.$$

Let us start with a

\subsubsection{ Multitime optimal control problem of \\maximizing a multiple integral functional}

{\it Find}
$$\max_{u(\cdot)}\,I(u(\cdot)) = \int_{\Omega_{t_0t_1}}L(t,x(t),u(x(t)))\, \om \leqno(27)$$
{\it subject to}
$$dx(t) = u^a_\al(x(t))X_a(x(t))dt^\al,\,\hbox{a.e.}\,\,\, t\in \Omega_{t_0t_1},\, x(t_0) = x_0,\, x(t_1)= x_1.\leqno(28)$$
It is supposed that $L: \Omega_{t_0t_1} \times A\times U \to R$ is a $C^2$ function and
$u^a: A \to R^{n-1}, \, a = 1,..., n-1$, $X^i_a: A \to R^{(n-1)n}, \, a = 1,..., n-1, i= 1,...,n$
are $C^2$ functions. Ingredients: $\om = dt^1\cdots dt^m$ is the volume element, $A$ is a bounded and closed subset of $R^n$,
which the $m$-sheet of controlled system is constrained to stay
for $t\in \Omega_{t_0t_1}$, and $x_0$ and $x_1$ are the initial and final states of the $m$-sheet $x(t)$ in the controlled system.
The set in which the control functions $u^a(t)$ takes their values in it, is called as $U$, which is a bounded and
closed subset of $R^{n-1}$.

Let us find the first order necessary conditions for an optimal pair $(x,u)$. Firstly, the {\it multitime infinitesimal deformation (Pfaff) system}
of the constraint (28) is (19), and the {\it multitime adjoint Pfaff system} is (20).

The control variables may be {\it open-loop} $u^a(t)$, depending directly on the multitime variable $t$, or {\it closed-loop}
(or {\it feedback}) $u^a(x(t))$, depending on the state $x(t)$.

To simplify, we accept an open-loop control. Introducing the $(m-1)$-forms
$$\omega_\la =\di \frac{\pa}{\pa t^\la}\rfloor\omega,$$
a {\it costate variable matrix or Lagrange multiplier matrix} $p=p^\al_i(t)\di\frac{\pa}{\pa t^\al}\otimes dx^i$ is identified to
the $(m-1)$-forms $p_i = p^\la_i(t) \omega_\la$.
We use the Lagrangian $m$-form
$${\cal L}(t,x(t),u(t),p(t)) = L(t,x(t),u(t))\om+p_i^\la(t)[u^a_\al(t)X^i_a(x(t))dt^\al - dx^i(t)]\wedge \omega_\la$$
and the {\it Hamiltonian $m$-form}
$${\cal H} = L(t,x(t),u(t))\om+p_i^\la(t)u^a_\al(t)X^i_a(x(t))dt^\al \wedge \omega_\la.$$

{\bf Theorem} ({\bf Multitime maximum principle})
{\it Suppose that the problem of maximizing the functional (27) constrained by (28) has
an interior optimal solution $\hat u(t)$, which determines the optimal evolution $x(t)$. Then there exists a
costate matrix $(p(t))$ such that
$$
dx\wedge \omega_\la = \frac{\partial {\mathcal H}}{\partial p^\la}, \leqno(29)
$$
the function $(p(t))$ is the unique solution of the following Pfaff system (adjoint system)
$$
dp^\la_j\wedge \omega_\la = - \frac{\pa {\cal H}}{\pa x^j},\,\,\, \delta_{\al\be}p^\al_i(t) n^\be(t) = 0\leqno(30)
$$
and the critical point conditions
$$
{\mathcal H}_{u^a}\left( t,x(t),u(t),p(t)\right) =0,\,a = \overline{1,n-1}\leqno(31)
$$
hold.}

{\bf Proof} We use the Lagrangian $m$-form ${\cal L}$.
The solutions of the foregoing problem are between the solutions of the free maximization problem of the functional
$$J(u(\cdot)) = \int_{\tilde\Omega}{\cal L}(t,x(t),u(t),p(t)),$$
where $\tilde \Omega = (\Omega_{t_0t_1}, x(\Omega_{t_0t_1}))\subset R^m_+\times R^n$.

Suppose that there exists a continuous control $\hat u(t)$ defined over the interval
$\Omega_{t_0t_1}$ with $\hat u(t)\in \hbox{Int}\,{U}$
which is an optimum point in the previous problem. Now consider a variation
$u(t,\epsilon)=\hat u(t) + \epsilon h(t)$, where $h$ is an arbitrary continuous vector function.
Since $\hat u(t)\in \hbox{Int}\,{\cal U}$ and a continuous function over a compact set
$\Omega_{t_0t_1}$ is bounded, there exists $\epsilon_h>0$ such
that $u(t,\epsilon)=\hat u(t) + \epsilon h(t)\in \hbox{Int}\,{U},\,\,\forall |\epsilon|<\epsilon_h$.
This $\epsilon$ is used in our variational arguments.

Define $x(t,\epsilon)$ as the $m$-sheet of the state variable corresponding to the control variable $u(t,\epsilon)$, i.e.,
$$dx^i(t;\epsilon) = u^a_\al(x(t;\epsilon))X^i_a(x(t;\epsilon))dt^\al,\, x(t;0) = x(t).$$
and $x(t_0,0)=x_0$, $x(t_1,0)=x_1$. For $|\epsilon|<\epsilon_h$, we define the function
$$J(\epsilon)= \int_{\tilde\Omega(\ep)}{\cal L}(t,x(t,\epsilon),u(t,\epsilon),p(t))$$
$$= \int_{\tilde\Omega(\ep)}L(t,x(t),u(t))\om+p_i^\la(t)[u^a_\al(t)X^i_a(x(t))dt^\al - dx^i(t)]\wedge \omega_\la.$$
Differentiating with respect to $\epsilon$, it follows
$$J^{\prime}(\epsilon)= \int_{\tilde\Omega(\ep)}\left(L_{x^j}(t,x(t,\epsilon),u(t,\epsilon)) + p_i^\al(t)u^a_\al(t,\epsilon) X^i_{ax^j}(x(t,\epsilon)\right)x^j_{\epsilon}(t,\epsilon)\om$$
$$ -  \int_{\tilde\Omega(\ep)} p_i^\la(t) dx^i_{\epsilon}(t,\epsilon)\wedge \omega_\la$$
$$ + \int_{\tilde\Omega(\ep)}\left({L}_{u^a_\la}(t,x(t,\epsilon),u(t,\epsilon) + p_i^\la(t) X^i_a(x(t,\epsilon))\right)h^a_\la(t) \om.$$
Evaluating at $\epsilon = 0$, we find
$$J^{\prime}(0)= \int_{\tilde\Omega(\ep)}\left(L_{x^j}(t,x(t),\hat u(t)) + p_i^\al(t)\hat{u}^a_\al(t) X^i_{ax^j}(x(t))\right)x^j_{\epsilon}(t,0)\om$$
$$ -  \int_{\tilde\Omega(\ep)} p_i^\la(t) dx^i_{\epsilon}(t,0)\wedge \omega_\la$$
$$ + \int_{\tilde\Omega(\ep)}\left({L}_{u^a_\la}(t,x(t),\hat u(t) + p_i^\la(t) X^i_a (x(t))\right)h^a_\la(t) \om.$$
where $x(t)$ is the $m$-sheet of the state variable corresponding to the optimal control $\hat u(t)$.

To evaluate the multiple integral
$$\int_{\tilde\Omega(\ep)} p^\la_i(t) dx^i_{\epsilon}(t,0)\wedge \omega_\la,$$
we integrate by parts, via the formula
$$d(p^\la_i x^i_{\epsilon} \omega_\la)=(x^i_\epsilon dp^\la_i + p^\la_i dx^i_\epsilon)\wedge\omega_\la,$$
obtaining
$$\di\int_{\tilde\Omega(\ep)}p^\la_idx^i_\epsilon\wedge \omega_\la=\di\int_{\tilde\Omega(\ep)}d(p^\la_ix^i_\epsilon\omega_\la)
- \di\int_{\tilde\Omega(\ep)} x^i_\epsilon dp^\la_i\wedge \omega_\la.$$
Now we apply the Stokes integral formula
$$\di\int_{\tilde\Omega(\ep)}d(p^\la_ix^i_\epsilon\omega_\la) = \di\int_{\pa{\tilde\Omega(\ep)}}\delta_{\al\be}p^\al_ix^i_\epsilon n^\be d\sigma,$$
where $(n^\be(t))$ is the unit normal vector to the boundary ${\pa{\tilde\Omega}}$.
Since the integral from the middle can be written
$$\di\int_{\tilde\Omega(\ep)}p^\la_idx^i_\epsilon\wedge \omega_\la=\di\int_{\pa{\tilde\Omega(\ep)}}\delta_{\al\be}p^\al_ix^i_\epsilon n^\be d\sigma
- \di\int_{\tilde\Omega} x^i_\epsilon dp^\la_i\wedge \omega_\la,$$
we find $J^{\prime}(0)$ as
$$J^{\prime}(0)= \int_{\tilde\Omega}\left(L_{x^j}(t,x(t),\hat u(t)) + p_i^\al(t)\hat{u}^a_\al(t) X^i_{ax^j}(x(t))\right)x^j_{\epsilon}(t,0)\om$$
$$- \di\int_{\pa{\tilde\Omega}}\delta_{\al\be}p^\al_ix^i_\epsilon(t,0) n^\be d\sigma + \di\int_{\tilde\Omega} dp^\la_i\wedge \omega_\la x^i_\epsilon(t,0)$$
$$ + \int_{\tilde\Omega}\left({L}_{u^a_\la}(t,x(t),\hat u(t) + p_i^\la(t) X^i_a (x(t))\right)h^a_\la(t) \om.$$
We select the costate function $p(t)$ as solution of the adjoint Pfaff equation (boundary value problem)
$$\left(L_{x^j}(t,x(t),\hat u(t)) + p_i^\al(t)\hat{u}^a_\al(t) X^i_{ax^j}(x(t))\right)\om + dp^\la_j\wedge \omega_\la = 0, \delta_{\al\be}p^\al_i(t) n^\be(t)|_{\pa\Omega} = 0.$$
On the other hand, we need $J^{\prime}(0)=0$ for all $h(t)=(h^a(t))$. Since
the variation $h$ is arbitrary, we get (critical point condition)
$${L}_{u^a_\la}(t,x(t),\hat u(t) + p_i^\la(t) X^i_a (x(t)) = 0.$$

Let us start with a

\subsubsection{Multitime optimal control problem of \\maximizing a curvilinear integral functional}

{\it Find}
{$$\max_{u(\cdot)}\,I(u(\cdot)) = \int_{\Gamma_{t_0t_1}}L_\al(t,x(t),u(x(t))) dt^\al \leqno(32)$$
{\it subject to}
$$dx(t) = u^a_\al(x(t))X_a(x(t))dt^\al\,\hbox{a.e.}\,\,\, t\in \Omega_{t_0t_1},\, x(t_0) = x_0,\, x(t_1)= x_1.\leqno(33)$$}
It is supposed that $L_\al: \Omega_{t_0t_1} \times A\times U \to R$ and $u^a: A \to R^{n-1}, \, a = 1,..., n-1$, $X^i_a: A \to R^{(n-1)n}, \, a = 1,..., n-1, i= 1,...,n$
are $C^2$ functions. Ingredients: ${\cal L} = L_\al(t,x(t),u(t)) dt^\al$ is an $1$-form, $A$ is a bounded and closed subset of $R^n$,
which the $m$-sheet of controlled system is constrained to stay
for $t\in \Omega_{t_0t_1}$, and $x_0$ and $x_1$ are the initial and final states of the $m$-sheet $x(t)$ in the controlled system.
The set in which the control functions $u^a$ takes their values in it, is called as $U$, which is a bounded and
closed subset of $R^{n-1}$.

Let us find the first order necessary conditions for an optimal pair $(x,u)$. Firstly, the {\it multitime infinitesimal deformation (Pfaff) system}
of the constraint (33) is (19), and the {\it multitime adjoint Pfaff system} is (*).

The control variables may be {\it open-loop} $u^a(t)$, depending directly on the multitime variable $t$, or {\it closed loop}
(or {\it feedback}) $u^a(x(t))$, depending on the state $x(t)$.

To simplify, we accept an open loop control. Introducing a {\it costate variable vector or
Lagrange multiplier} $p = (p_i(t))$, we build a Lagrangian $1$-form
$${\cal L}(t,x(t),u(t),p(t)) = L_\al(t,x(t),u(t)) dt^\al + p_i(t)[u^a_\al(t)X^i_a(x(t))dt^\al - dx^i(t)].$$
and a {\it Hamiltonian $1$-form}
$${\cal H} = L_\al(t,x(t),u(t)) dt^\al + p_i(t)u^a_\al(t)X^i_a(x(t))dt^\al.$$

{\bf Theorem} ({\bf Multitime maximum principle})
{\it Suppose that the problem of maximizing the functional (32) constrained by (33) has
an interior optimal solution $\hat u(t)$, which determines the optimal evolution $x(t)$. Then there exists a
costate function $(p(t))$ such that
$$
dx^i = \frac{\partial {\mathcal H}}{\partial p_i}, \leqno(34)
$$
the function $(p(t))$ is the unique solution of the following Pfaff system (adjoint system)
$$
dp_i = - \frac{\pa {\cal H}}{\pa x^i}\leqno(35)
$$
and the critical point conditions
$$
{\mathcal H}_{u^a_\al}\left( t,x(t),u(t),p(t)\right) =0,\,a = \overline{1,n-1}, \al = 1,...,m\leqno(36)
$$
hold.}

{\bf Proof} We use the Lagrangian $1$-form ${\cal L}$.
The solutions of the foregoing problem are between the solutions of the free maximization problem of the curvilinear integral functional
$$J(u(\cdot)) = \int_{\tilde\Gamma}{\cal L}(t,x(t),u(t),p(t)),$$
where $\tilde \Gamma = (\Gamma_{t_0t_1}, x(\Gamma_{t_0t_1}))\subset R^m_+\times R^n$.


Suppose that there exists a continuous control $\hat u(t)$ defined over the interval $\Omega_{t_0t_1}$ with $\hat u(t)\in \hbox{Int}\,{U}$
which is an optimum point in the previous problem. Now consider a variation
$u(t,\epsilon)=\hat u(t) + \epsilon h(t)$, where $h$ is an arbitrary continuous vector function.
Since $\hat u(t)\in \hbox{Int}\,{\cal U}$ and a continuous function over a compact set
$\Omega_{t_0t_1}$ is bounded, there exists $\epsilon_h>0$ such
that $u(t,\epsilon)=\hat u(t) + \epsilon h(t)\in \hbox{Int}\,{U},\,\,\forall |\epsilon|<\epsilon_h$.
This $\epsilon$ is used in our variational arguments.

Define $x(t,\epsilon)$ as the $m$-sheet of the state variable corresponding to the control variable $u(t,\epsilon)$, i.e.,
$$dx^i(t;\epsilon) = u^a_\al(x(t;\epsilon))X^i_a(x(t;\epsilon))dt^\al,\, x(t;0) = x(t)$$
and $x(t_0,0)=x_0$, $x(t_1,0)=x_1$. For $|\epsilon|<\epsilon_h$, we define the function
$$J(\epsilon)= \int_{\tilde\Gamma(\ep)}{\cal L}(t,x(t,\epsilon),u(t,\epsilon),p(t))$$
$$= \int_{\tilde\Gamma(\ep)}L_\al(t,x(t,\epsilon),u(t,\epsilon)) dt^\al + p_i(t)[u^a_\al(t,\epsilon)X^i_a(x(t,\epsilon))dt^\al - dx^i(t,\epsilon)].$$
Differentiating with respect to $\epsilon$, it follows
$$J^{\prime}(\epsilon)= \int_{\tilde\Gamma(\ep)}\left({L}_{\al x^j}(t,x(t,\epsilon),u(t,\epsilon)) dt^\al+ p_i(t)u^a_\al(t,\epsilon)X^i_{ax^j}(x(t,\epsilon))dt^\al\right)x^j_{\epsilon}(t,\epsilon)$$
$$ -  \int_{\tilde\Gamma(\ep)} p_i(t) dx^i_{\epsilon}(t,\epsilon)$$
$$ + \int_{\tilde\Gamma(\ep)}\left({L}_{\be u^a_\al}(t,x(t,\epsilon),u(t,\epsilon))dt^\be + p_i(t)X^i_a(x(t,\epsilon))dt^\al\right) h^a_\al(t).$$
Evaluating at $\epsilon = 0$, we find $\tilde\Gamma(0) = \tilde\Gamma$ and
$$J^{\prime}(0)= \int_{\tilde\Gamma}\left({L}_{\al x^j}(t,x(t),\hat u(t)) dt^\al+ p_i(t)\hat {u}^a_\al(t)X^i_{ax^j}(x(t))dt^\al\right)x^j_{\epsilon}(t,0)$$
$$ -  \int_{\tilde\Gamma} p_i(t) dx^i_{\epsilon}(t,0)$$
$$ + \int_{\tilde\Gamma}\left({L}_{\be u^a_\al}(t,x(t),\hat u(t))dt^\be + p_i(t)X^i_a(x(t))dt^\al\right) h^a_\al(t).$$
where $x(t)$ is the $m$-sheet of the state variable corresponding to the optimal control $\hat u(t)$.

To evaluate the curvilinear integral
$$\int_{\tilde\Gamma} p_i(t) dx^i_{\epsilon}(t,0),$$
we integrate by parts, via the formula
$$d(p_i x^i_{\epsilon})= x^i_\epsilon dp_i + p_i dx^i_\epsilon,$$
obtaining
$$\di\int_{\tilde\Gamma}p_idx^i_\epsilon = (p_i(t) x^i_{\epsilon}(t,0))|_{t_0}^{t_1} - \di\int_{\tilde\Gamma}(dp_i) x^i_\epsilon.$$
We find

$$J^{\prime}(0)= \int_{\tilde\Gamma}\left({L}_{\al x^j}(t,x(t),\hat u(t)) dt^\al+ p_i(t)\hat {u}^a_\al(t)X^i_{ax^j}(x(t))dt^\al+ dp_j\right)x^j_{\epsilon}(t,0)$$
$$ -  (p_i(t) x^i_{\epsilon}(t,0))|_{t_0}^{t_1} $$
$$ + \int_{\tilde\Gamma}\left({L}_{\be u^a_\al}(t,x(t),\hat u(t))dt^\be + p_i(t)X^i_a(x(t))dt^\al\right) h^a_\al(t).$$
We select the costate function $p(t)$ as solution of the adjoint Pfaff equation (boundary value problem)
$${L}_{\al x^j}(t,x(t),\hat u(t)) dt^\al+ p_i(t)\hat {u}^a_\al(t)X^i_{ax^j}(x(t))dt^\al+ dp_j = 0.$$
On the other hand, we need $J^{\prime}(0)=0$ for all $h(t)=(h^a_\al(t))$. Since
the variation $h$ is arbitrary, we get (critical point condition)
$${L}_{\be u^a_\al}(t,x(t),\hat u(t))dt^\be + p_i(t)X^i_a(x(t))dt^\al = 0.$$

{\bf Example: Nonholonomic control of torsion of a cylinder or prism}

Suppose the torsion of a cylinder or prism is described by the controlled Pfaff equation
$$dz = (y + u(x,y))dx + (-x + v(x,y))dy,\,\,z(0,0) = 0, z(x_0,y_0) = z_0,$$
where the control $(u,v)$ is not subject to constraints. If the complete integrability condition $\di{\pa v\ov \pa x}-\di{\pa u\ov \pa y} = 2$
is verified identically, then we have a holonomic evolution. Otherwise, we have a nonholonomic evolution. Using
the controlled Pfaff equation as constraint, we want to minimize the functional
$$\int_\Gamma (z(x,y)+u(x,y)^2)c_1dx+(z(x,y)+v(x,y)^2) c_2dy,$$
where $\Gamma$ is a $C^1$ curve joining the points $(0,0)$ and $(x_0,y_0)$, and $c_1, c_2$ are constants.
The minimization of the previous integral is equivalent to the maximization of the cost functional
$$P(u(\cdot),v(\cdot))= - \int_\Gamma (z(x,y)+u(x,y)^2)c_1dx+(z(x,y)+v(x,y)^2) c_2dy$$
subject to controlled Pfaff equation.

Let us find the optimal manifold (surface or curve) of evolution, using the two-variable maximum principle theory.
For that we introduce the 1-forms:
$$\begin{array}{c}\omega = (y+u)dx+(-x+v)dy - dz\,\, (\hbox{evolution 1-form})\\ \
{\omega}^0 = - (z+u^2)c_1 dx - (z+v^2)c_2 dy\,\, (\hbox{running cost 1-form})\\ \
\eta = p_0 \omega^0 + p\, \omega\,\,(\hbox{control 1-form}).
\end{array}$$
Taking $p_0=1$, we obtain
$$\eta = -(z+u^2)c_1 dx - (z+v^2)c_2 dy) + p((y+u)dx+(-x+v)dy - dz).$$
The adjoint equation $dp = -\di{\pa \ov \pa z}\eta =c_1 dx+c_2 dy,\,\,p(0,0)=0$ has the solution $p(x,y)=c_1 x+c_2 y$.
The maximization condition
$$H_1(x,y,z(x,y),p(x,y),u(x,y)) =\max_{u,v}\{-(z(x,y)+u^2)c_1 + p(x,y)(y+u)\}$$
$$H_2(x,y,z(x,y),p(x,y),v(x,y)) =\max_{u,v}\{-(z(x,y)+v^2)c_1 + p(x,y)(y+u)\}$$
gives the optimal law
$$u(x,y) = {p(x,y)\ov 2c_1},\,v(x,y) = {p(x,y)\ov 2c_2}.$$
Replacing in the evolution Pfaff equation we obtain
$$dz = \left(y + {c_1 x+c_2 y\ov 2c_1}\right)dx + \left(-x + {c_1 x+c_2 y\ov 2c_2}\right)dy.$$

1) If the complete integrability condition
$${\pa\ov \pa y}\left(y + {c_1 x+c_2 y\ov 2c_1}\right)={\pa\ov \pa x}\left(-x + {c_1 x+c_2 y\ov 2c_2}\right),$$
i.e., $4 c_1c_2 + c_2^2 -c_1^2 = 0$ is satisfied, then the evolution surface is
$$z = xy + \frac{x^2}{4} +\frac{y^2}{4}+\frac{c_2}{2c_1}xy.$$

2) If $4 c_1c_2 + c_2^2 -c_1^2 \neq 0$, then the Pfaff evolution equation admits only solutions
which are curves (nonholonomic surface in $R^3$):
$x = x(t), y = y (t), z = z(t), t\in I,$ with
$$\begin{array}{c}x = x(t), y = y(t) (\hbox{given arbitrary})\\ \
\di{dz\ov dt} =y(t){x}'(t) - x(t) {y}'(t) + (c_1 x(t) +c_2 y(t)){\left(\frac{{x}'(t)}{2c_1} +\frac{{y}'(t)}{2c_2}\right)}.
\end{array}$$
In this case, for determining $z = z(t)$, we must take $x = x(t), y = y (t),\, t\in I$ as a parametrization of the curve $\Gamma$
from the cost functional. In fact, the problem is reduced to optimization of a simple integral constrained by
a differential equation.


\section{Curvilinear integral functionals \\depending on curves}

Let $\Omega_{t_0t_1}\subset R^m$, let ${\Gamma_{t_0t_1}}\subset \Omega_{t_0t_1}$ be a $C^1$ curve and
$$J({\Gamma_{t_0t_1}}; x(\cdot)) = \int_{\Gamma_{t_0t_1}}\,L_\al(t,x(t),x_\ga(t)) dt^\al$$
be a curvilinear integral functional depending on the curve ${\Gamma_{t_0t_1}}$. Consider a variation
${\Gamma_{t_0t_1}(\ep)}: t = (t^\al(\tau,\ep))$ of the curve ${\Gamma_{t_0t_1}}: t = (t^\al(\tau))$, with
the same endpoints. Suppose $L_\al(t,x(t),x_\ga(t)) dt^\al$ is stationary with respect to $\ep$.
Then $$J(\ep) = \int_{\Gamma_{t_0t_1}(\ep)}\,L_\al(t,x(t),x_\ga(t)) dt^\al.$$

The closed curve $C = \Gamma_{t_0t_1}(\ep)\cup {\Gamma_{t_1t_0}(0)}$ is the boundary of a surface $S$.
We evaluate $J(\ep) - J(0)$, using Stokes formula,
$$J(\ep) - J(0) = \int_{\Gamma_{t_0t_1}(\ep)}\,L_\al(t,x(t),x_\ga(t)) dt^\al - \int_{\Gamma_{t_0t_1}(0)}\,L_\al(t,x(t),x_\ga(t)) dt^\al$$
$$ = \int_{\Gamma_{t_0t_1}(\ep)}\,L_\al(t,x(t),x_\ga(t)) dt^\al + \int_{\Gamma_{t_1t_0}(0)}\,L_\al(t,x(t),x_\ga(t)) dt^\al$$
$$ = \int_C  \,L_\al(t,x(t),x_\ga(t)) dt^\al = \int_S \, d(L_\al(t,x(t),x_\ga(t)) dt^\al) = \int_S \, D_\be L_\al dt^\be\wedge dt^\al$$
$$= \frac{1}{2}\int_S (D_\al L_\be - D_\be L_\al)dt^\al\wedge dt^\be.$$
Now we use the variation vector field $\frac{\pa t^\al}{\pa\ep}(\tau, \ep)|_{\ep = 0} = \xi^\al(\tau)$. Replacing $dt^\al =\ep\, \xi^\al$, the surface integral
is transformed to a curvilinear integral
$$= \frac{\ep}{2}\int_{\Gamma_{t_0t_1}(0)}\, (D_\al L_\be - D_\be L_\al)(\xi^\al dt^\be - \xi^\be dt^\al) =
\ep\,\int_{\Gamma_{t_0t_1}(0)}\, (D_\al L_\be - D_\be L_\al)\xi^\al dt^\be.$$
It follows
$$J^{\prime}(0) = \int_{\Gamma_{t_0t_1}(0)}\, (D_\al L_\be - D_\be L_\al)\xi^\al dt^\be.$$

Suppose ${\Gamma_{t_0t_1}}$ is a critical point of the functional, hence $J^{\prime}(0) = 0,  \,\forall \xi$.
Consequently
$$(D_\al L_\be - D_\be L_\al)(t(\tau))\, \frac{dt^\be}{d\tau}(\tau) = 0.$$
If the curvilinear integral is path independent, then this relation is identically satisfied.
If the curvilinear integral is path dependent, then the discussion depends on $m$
since $a_{\al\be} = D_\al L_\be - D_\be L_\al$ is an anti-symmetric matrix, and consequently its determinant $d = \det(a_{\al\be})$
is either $0$, for $m$ odd, or $\geq 0$, for $m$ even. For $m$-odd we have solutions, i.e., critical curves;
for $m$-even, we have either no solution for $d > 0$ or solutions for $d = 0$.
Since the differential system is of order one, the curve solution is determined only by a single condition
(the general bilocal problems have no solution). The extremum problems have sense only if we add supplimentary conditions
(an initial condition + an isoperimetric condition).

{\bf Variant} Let $\Omega_{t_0t_1}\subset R^m$, let ${\Gamma_{t_0t_1}}\subset \Omega_{t_0t_1}$ be a $C^1$ curve and
$$J({\Gamma_{t_0t_1}}; x(\cdot)) = \int_{\Gamma_{t_0t_1}}\,L_\al(t,x(t),x_\ga(t)) dt^\al$$
be a curvilinear integral functional depending on the curve ${\Gamma_{t_0t_1}}$. Consider a variation
${\Gamma_{t_0t_1}(\ep)}: t = (t^\al(\tau,\ep))$ of the curve ${\Gamma_{t_0t_1}}: t = (t^\al(\tau))$, with
the same endpoints. Denote $M_\al(t) = L_\al(t,x(t),x_\ga(t))$.

Then
$$J(\ep) = \int_{\Gamma_{t_0t_1}(\ep)}\,M_\al(t(\tau,\ep)) dt^\al(\tau, \ep).$$
To compute $J^{\prime}(0)$, we use the variation vector field $\frac{\pa t^\al}{\pa \ep}(\tau, \ep)|_{\ep = 0} = \xi^\al(\tau)$.
From
$$J^{\prime}(\ep) = \int_{\Gamma_{t_0t_1}(\ep)}\frac{\pa M_\al}{\pa t^\be}(t(\tau, \ep))\frac{\pa t^\be}{\pa\ep}(\tau,\ep)dt^\al(\tau, \ep) +
M_\al(t(\tau, \ep))\,d \,\,\frac{\pa t^\al}{\pa \ep}(\tau, \ep),$$
we obtain
$$J^{\prime}(0) = \int_{\Gamma_{t_0t_1}(0)}\frac{\pa M_\al}{\pa t^\be}(t(\tau))\xi^\be(\tau)dt^\al(\tau) + M_\be(t(\tau))d \xi^\be(\tau).$$
Integrating by parts, we find
$$J^{\prime}(0) =  M_\be(t(\tau))\xi^\be(\tau)|^{\tau_1}_{\tau_0} +\int_{\Gamma_{t_0t_1}(0)}\left(\frac{\pa M_\al}{\pa t^\be}- \frac{\pa M_\be}{\pa t^\al}\right)(t(\tau))\xi^\be(\tau)dt^\al(\tau).$$

{\bf Remark} The variation of the function $x(t)$ has nothing to do with the variation of the curve.

\section{Optimization of mechanical work on \\Riemannian manifolds}

Let $(M,g)$ be a Riemannian manifold and $X$ a $C^{2}$ vector field on M.
Let $x=(x^1,...,x^n)$ denote the local coordinates relative to a fixed local map $(V,h)$.
Since $h:V\rightarrow R^n$ is a diffeomorphism, we denote by $\Omega_{x_{0}x_{1}}$
a subset of $V$ diffeomorphic through $h$ with the
hyper-parallelepiped in $R^n$ having $h(x_{0})$ and $h(x_1)$ as diagonal points.

Let ${\Gamma_{x_0x_1}}: x^i = x^i(t), t\in [t_0,t_1]$ be an arbitrary $C^1$ curve on $M$ which joins the points $x(t_0) = x_0, x(t_1) = x_1.$
The functional
$$J(\Gamma_{x_0x_1}) = \int_{\Gamma_{x_0x_1}}\, g_{ij}(x)X^i(x) dx^j$$
is generated by the {\it mechanical work} produced by the force $\om_j = g_{ij}(x)X^i(x)$ along the curve $\Gamma_{x_0x_1}$.

Let $X$ be a nowhere zero vector field. $X$ is called a {\it geodesic vector field} iff $\nabla_X X = 0$. Thus $X$ is geodesic iff
each of its integral curves is a geodesic.

{\bf Theorem} {\it If $X$ is a unit geodesic vector field and $\gamma_{x_0x_1}$ is a field line,
then the curve $\gamma_{x_0x_1}$ is a maximum point of the functional $J(\Gamma_{x_0x_1})$
and the maximum value is the length of $\gamma_{x_0x_1}$.}

{\bf Proof} Let us find
$$\max_{\Gamma_{x_0x_1}} J(\Gamma_{x_0x_1}) = \int_{\Gamma_{x_0x_1}}\, g_{ij}(x)X^i(x) dx^j,$$
where
$$g_{ij}(x)X^i(x)X^j(x) = 1.$$

The critical point condition, with respect to the curve ${\Gamma_{x_0x_1}}$, is
$$\left(g_{ij} \nabla_kX^i - g_{ik}\nabla_j X^i\right)(x(t)) \frac{dx^j}{dt}(t) = 0.$$
It is identically satisfied, because ${\gamma_{x_0x_1}}$ is a field line,
the geodesic condition implies $\nabla_X X^i = 0$ and the the condition of unit vector field gives
$g_{ij}X^j \nabla_k X^i = 0.$

On the other hand, the inequality
$$|g_{ij}(x)X^i(x) dx^j| \leq ||X|| ds,\, ds = ||dx|| = \sqrt{g_{ij}dx^i dx^j}$$
becomes an equality if $dx^j = X^j(x(t))dt$, i.e., $x(t)$ is a field line of $X(x)$.
Under the condition $||X|| =1$, the maximum value of the foregoing functional is the length of ${\gamma_{x_0x_1}}$.

\section{Bang-bang control on distributions}

The same distribution $D$ can be described in terms of vector fields,
$$D = \hbox{span}\{X_a(x)|\, a_i(x)X^i_a = 0,\, a = 1,...,n-1\}.\leqno (22)$$

Bang-bang control is an optimal or suboptimal piecewise constant control whose values
are defined by bounds imposed on the amplitude of control components.
The control changes its values according to the switching function which may be found using the maximum principle.
The discontinuity of the bang-bang control leads to discontinuity of a value function for the considered optimal control problem.
Typical problems with bang-bang optimal control include time and terminal cost optimal control for linear control systems.
Bang-bang optimization offers a direct explanation for an otherwise perplexing observation and indicates
that evolution is operating according to principles that every engineer knows.
The balck hole applications covered in this Section refer to the controllability of the ODE or PDE system by bang-bang controls.

\subsection{Single-time bang-bang optimal control}

Let $x(t),\,\,t\in I = [0, \tau] \subset R$, be an integral curve of the distribution $D$.
Any curve in the distribution $\Delta = \hbox{span}\{X_a, a = 1,...,n-1\}$ is a solution of the controlled ODE system
$$\dot x(t) = u^a (t)X_a (x(t)),\,\,u(t) = (u^a(t)), \,\, t\in [0, \tau],\eqno(ODE)$$
called {\it driftless control system}.

{\bf (1) Time minimum problem} Let $U = [-1, 1]^{n-1} \subset R^{n-1}$ be the control set. Giving the starting point $x_0 \in
R^n$, find an optimal control $u^*(\cdot)$ such that
$$
I(u^* (\cdot)) = \min_{u(\cdot)} \int_{0}^{\tau}\, dt,
$$
using (ODE) evolution as constraint. Since
$\tau^{*} = I(u^*(\cdot))$, the optimal point  $\tau^*$ ensures the minimum time
to steer to the origin. This time optimum problem is equivalent to a controllability one.

{\bf Solution} To prove the existence of a bang-bang control, we use the single-time Pontryaguin Maximum Principle.
The Hamiltonian $H(x,p,u) = - 1+ p_i X^i_a(x) u^a$ gives the adjoint ODE system $\dot p_j(t) = - p_i(t) \frac{\pa X^i_a}{\pa x^j}(x(t)) u^a(t)$.
The extremum of the linear function $u\to H$ exists since each control variable belong to
the interval $[- 1, 1]$; for optimum, the control must be at a vertex of $\pa U$ (see, linear optimization, simplex method). If
$Q_a(t) = p_i(t) X^i_a(x(t))$, then the optimal control $u^{*a}$ must be the function (bang-bang control)
$$
u^{*a} = -\, \hbox{sign}\, (Q_{a}(t)) =\left\{\begin{array}{cc} \hspace{-2.8cm}1 \hspace{2.9cm}\hbox{for}\,\,\, Q_a(t) < 0\\ \
\hspace{-0.2cm}\hbox{undetermined\ \hspace{0.8cm} for}\,\, Q_a (t)= 0 \,(\hbox{singular \,control}) \\ \
\hspace{-3.cm} - 1\hspace{2.7cm} \hbox{for}\,\, Q_a(t) > 0.\end{array}\right.
$$
Suppose the Lebesgue measure of each set $\{t \in [0, \tau] : Q_{a}(t) = 0\}$ vanishes.
Then the singular control is ruled out and the remaining possibilities are bang-bang controls.
This optimal control is discontinuous since each component
jumps from a minimum to a maximum and vice versa in response to each change
in the sign of each $Q_a(t)$. The functions $Q_a(t)$ are called {\it switching functions}.

{\bf (2) Optimal terminal value} Let $U = [-1, 1]^{n-1} \subset R^{n-1}$ be the control set.
Suppose we have to

{\it Minimize the terminal cost functional
$$Q(u(\cdot)) = x^n(t_1)$$
subject to the driftless control system
$$
\dot x(t) = u^a(t)X_a (x(t)),\,\,u(t)\in {\cal U}, \,\, t\in [t_0, t_1];\,\, x(t_0)=x_0.
$$}

{\bf Solution} Since the control Hamiltonian $H(x,p,u) = p_i X^i_a(x) u^a$ is linear in the control, the optimal control is a bang-bang.
Automatically we find the optimal costate function and the optimal evolution.


\subsection{Multitime bang-bang optimal control}

Let $\Omega_{0\tau}$ be the hyperparallelipiped determined by two opposite diagonal points $0 =
(0,..., 0)$ and $\tau = (\tau^1, ...,\tau^m)$ in $R^m_+$, endowed with the product order.
Let $x(t),\,\,t\in \Omega_{0\tau} \subset R^m_+$, be an integral $m$-sheet of the distribution $D$, i.e.,
a solution of a multitime piecewise completely integrable PDE system
$$\frac{\pa x}{\pa t^\al}(t) = u^a_\al (t)X_a (x(t)),\,t = (t^\al)\in \Omega_{0\tau}, \, a = \overline{1,n-1}, \,\al = 1,...,m.\eqno(PDE)$$
This sort of controlled PDE is called a {\it driftless control system}. Of course, the piecewise complete integrability conditions
$$\left(\frac{\pa u^a_\al}{\pa t^\be} - \frac{\pa u^a_\be}{\pa t^\al}\right) X_a = u^a_\al u^b_\be [X_a, X_b]\eqno(CIC)$$
restrict the controls, excepting the case when they are identically satisfied.

To show that the driftless control system is multitime controllable, by bang-bang controls (see also [10]), we use the next
{\it multitime minimum problems}

{\bf Case of multiple integral functional}
Let $U = [-1, 1]^{m(n-1)} \subset R^{m(n-1)}$ be the control set. Giving the starting point $x_0 \in
R^n$, find an optimal control $u^*(\cdot)$ such that
$$
I(u^* (\cdot)) = \min_{u(\cdot)} \int_{\Omega_{0\tau}}\, dt^1... dt^m,
$$
using a completely integrable two-time evolution (PDE) as constraint and supposing that (CIC) are identically satisfied.
Since $\tau^{*1}\cdots \tau^{*2} = I(u^*(\cdot))$, the optimal point  $\tau^* = (\tau^{*1},..., \tau^{*m})$ ensures the minimum multitime
"volume" to steer to the origin. This two-time optimum problem consists in devising a control such that to transfer a given initial state
to a specified target (controllability problem).

{\bf Solution} We apply the multitime maximum principle which proves the existence of a bang-bang control.
The Hamiltonian $H(x,p,u) = - 1+ p_i^\al X^i_a(x) u^a_\al$ gives the adjoint PDE system
$\frac{\pa p^\al_j}{\pa t^\al}(t) = - p_i^\al(t) \frac{\pa X^i_a}{\pa x^j}(x(t)) u^a_\al(t)$.
The extremum of the linear function $u\to H$ exists since the set $U$ is compact; for optimum, the control vectors
$u_\al = (u^1_\al, ..., u^{n-1}_\al)$ must be vertices of $\pa U$. If
$Q^\al_a(t) = p^\al_i(t) X^i_a(x(t))$ are the switching functions, then each optimal control $u^{*a}_\al$ is of the form
$$
u^{*a}_\al = - \,\hbox{sign}\, (Q^\al_{a}(t)) =\left\{\begin{array}{cc}\hspace{0.6cm} 1 \hspace{1.9cm}\hbox{for}\,\,\, Q^\al_a(t) < 0: \,\hbox{bang-bang \,control}\\ \
\hspace{-0.6cm}\hbox{undetermined} \hspace{0.2cm} \hbox{for}\,\, Q^\al_a (t)= 0: \,\hbox{singular \,control} \\ \
\hspace{0.3cm}- 1\hspace{1.9cm} \hbox{for}\,\, Q^\al_a(t) > 0: \,\hbox{bang-bang \,control}.\end{array}\right.
$$
Suppose the Lebesgue measure of each set $\{t \in \Omega_{0\tau} : Q^\al_{a}(t) = 0\}$ vanishes.
Then the singular control is ruled out and the remaining possibilities are bang-bang controls.
This optimal control is discontinuous since each component
jumps from a minimum to a maximum and vice versa in response to each change
in the sign of each $Q^\al_a(t)$. The piecewise complete integrability identities keep only
the control vectors (vertices of $\pa U$) $u_\al$ which satisfy $u_\al = \pm\, u_1$.
Each optimal $m$-sheet $x(t)$ is a {\it soliton} solution.

{\bf Case of curvilinear integral functional}

{\bf Optimal terminal value} Let $U = [-1, 1]^{m(n-1)} \subset R^{m(n-1)}$ be the control set.
Suppose we have to

{\it Minimize the terminal cost functional
$$Q(u(\cdot)) = x^n(\tau)$$
subject to the driftless control system}
$$\frac{\pa x}{\pa t^\al}(t) = u^a_\al (t)X_a (x(t)),\,t = (t^\al)\in \Omega_{0\tau},\, x(0) = x_0, \, a = \overline{1,n-1}, \,\al = 1,...,m.\eqno(PDE)$$

{\bf Solution} Since the control Hamiltonian $H(x,p,u) = p_i X^i_a(x) u^a$ is linear in the control, the optimal control is a bang-bang.
Automatically we find the optimal costate function and the optimal evolution.

\section{Optimal control problems on \\Tzitzeica surfaces}

Let $R^2_+$ be endowed with the product order. Let $\Omega\subset R^2_+$ be the bi-dimensional interval
determined by the opposite diagonal points $(0,0)$ and $(u^1, v^1)$.

Problem: {\it find
$$\max_h \int_{\Omega}L(u,v,\vec r(u,v), h(u,v)) du dv$$
constrained by (non-ruled Tzitzeica surfaces)
$$\vec r_{uu} = \frac{h_u}{h}\vec r_u + \frac{1}{h}\vec r_v,\,\,\vec r_{uv} = h \vec r,\,\, \vec r_{vv} = \frac{1}{h}\vec r_u +\frac{h_v}{h}\vec r_v,\,\,
(\ln h)_{uv} = h - \frac{1}{h^2},$$
$$\vec r(0,0) = \vec r_0,\, \vec r(u^1,v^1) = \vec r_1,$$
where $h(u,v)$ is the control.}

To solve the problem we use
$${\cal L} = L + \langle \vec a, \frac{h_u}{h}\vec r_u + \frac{1}{h}\vec r_v - \vec r_{uu}\rangle + \langle \vec b, h \vec r - \vec r_{uv}\rangle$$
$$+ \langle \vec c, \frac{1}{h}\vec r_u +\frac{h_v}{h}\vec r_v - \vec r_{vv}\rangle + d\left(h - \frac{1}{h^2} - (\ln h)_{uv}\right).$$

Explicitely, find
$$\max_h \int_{\Omega}(h^2(u,v) + ||\vec r(u,v)||^2)du dv.$$

\section{Phytoplankton growth model }

{\bf Open problem} {\it Transform the next ODE systems into Pfaff systems and study
their stochastic perturbations}.

Alessandro Abate, Ashish Tiwari, Shankar Sastry, Box Invariance for biologically-inspired dynamical systems

(i) O. Bernard and J.-L. Gouze, "Global qualitative description of a class
of nonlinear dynamical systems," Artificial Intelligence, vol. 136, pp.
29-59, 2002:

Consider the following Phytoplankton Growth Model:
$$\dot x^1 = 1- x^1- \frac{1}{4}\,x^1x^2,\,\, \dot x^2 = 2x^2x^3 - x^2,\,\, \dot x^3 = \frac{1}{4}\,x^1 - 2{x^3}^2,$$
where $x^1$ denotes the substrate, $x^2$ the phytoplankton
biomass, and $x^3$ the intracellular nutrient per biomass.

(ii) A. Julius, A. Halasz, V. Kumar, and G. Pappas, "Controlling biological
systems: the lactose regulation system of Escherichia Coli," in
American Control Conference 2007:

The dynamics of tetracycline antibiotic in a
bacteria which develops resistance to this drug (by turning on
genes $TetA$ and $TetR$) can be described by the following hybrid
system:
$$\dot x^1 = f - \frac{1}{3}\,x^1 x^3 + \frac{1}{8000}\, x^2,\,\,\dot x^2 = \frac{3}{200}\,u_0 - \frac{7}{2}\,x^3x^4$$
$$\dot x^3 = \frac{1}{3}\,x^1x^3 - \frac{1}{2500}\,x^2,\,\,\dot x^4 = f - \frac{11}{40000}\,x^4,$$
where $f = \frac{1}{2000}$ if $TetR > \frac{1}{50000}$ and $f = \frac{1}{40}$
otherwise ($f$ is the transcription rate of genes, which are inhibited by $TetR$),
and $x^1$, $x^2$, $x^3$, $x^4$ are the cytoplasmic concentrations
of $TetR$ protein, the $TetR-Tc$ complex, Tetracycline,
and $TetA$ protein, and $u_0$ is the extracellular concentration
of Tetracycline.

{\bf Acknowledgments}

Partially supported by University Politehnica of Bucharest, by UNESCO Chair in Geodynamics, "Sabba S. \c Stef\u anescu"
Institute of Geodynamics, Romanian Academy and by Academy of Romanian Scientists.

Prof. Dr. Constantin Udri\c ste, University Politehnica of Bucharest, Faculty of Applied Sciences, Departament of Mathematics-Informatics,
Splaiul Independentei 313, 060042 Bucharest, Romania, \\E-mail: udriste@mathem.pub.ro, anet.udri@yahoo.com

\end{document}